\documentclass[12pt]{amsart}
\usepackage{latexsym, amssymb, amsmath, amscd}
 \usepackage{eucal}

    \newtheorem{rema}{Remark}[section]
    \newtheorem{propo}[rema]{Proposition}

   \newtheorem{theo}[rema]{Theorem}
   \newtheorem{def-theo}[rema]{Definition-Theorem}
 
   \newtheorem{defi}[rema]{Definition}
    \newtheorem{lemma}[rema]{Lemma}
    \newtheorem{corol}[rema]{Corollary}
     
  \newtheorem{rmk}[rema]{Remark}

	\newcommand{\nno}{\nonumber}

	\newcommand{\p}{\partial}

 \newcommand{\res}{\mbox{\rm Res}}

 \newcommand{\pf}{{\it Proof:}\hspace{2ex}}
 
 \newcommand{\epfv}{\hspace{1em}$\Box$\vspace{1em}}

\newcommand{\bC}{{\mathbb C}}

\newcommand{\bQ}{{\mathbb Q}}
\newcommand{\bN}{{\mathbb N}}
\newcommand{\bT}{{\mathbb T}}

\newcommand{\cT}{{\mathcal T}}

\newcommand{\cS}{{\mathcal S}}

\newcommand{\cD}{{\mathcal D}}
\newcommand{\cH}{{\mathcal H}}

\newcommand{\cNcs}{{${\mathcal N}$CS} }
\newcommand{\cNsf}{{{\mathcal N}Sym}}
\newcommand{\cQf}{{{\mathcal Q}Sym}}
\newcommand{\Oft}{ \Omega_{F_t}  }

\newcommand{\cSft}{ {\mathcal S}_{F_t} }

\newcommand{\cDz}{{\mathcal D \langle  z \rangle}}
\newcommand{\cDaz}{{{\mathcal D}^{[\alpha]} \langle  z \rangle}}

\newcommand{\cDzz}{{\mathcal D \langle \langle z \rangle\rangle}}
\newcommand{\cDkzz}{{\mathcal D_K \langle \langle z \rangle\rangle}}
\newcommand{\cDtzz}{{\mathcal D_t \langle \langle z \rangle\rangle}}
\newcommand{\cDrzz}{{\mathcal D er\langle \langle z \rangle\rangle}}
\newcommand{\cDrtzz}{{\mathcal D er_t \langle \langle z \rangle\rangle}}
\newcommand{\cDrkzz}{{\mathcal D er_K \langle \langle z \rangle\rangle}}

\newcommand{\cDrazz}{{\cD er^{[\alpha]}\langle \langle z \rangle \rangle}} 
 
\newcommand{\cDazz}{{\cD^{[\alpha]}\langle \langle z \rangle \rangle}}

\newcommand{\cDrtazz}{{\cD er^{[\alpha]}_t\langle \langle z \rangle \rangle}} 
 
\newcommand{\cDtazz}{{\cD^{[\alpha]}_t\langle \langle z \rangle \rangle}}

\newcommand{\ataz}{{\mathbb A_t^{[\alpha]}\langle \langle z\rangle\rangle}}
\newcommand{\gtaz}{{\mathbb G_t^{[\alpha]}\langle \langle z\rangle\rangle}}

\newcommand{\klam}{{K\langle \Lambda \rangle}}
\newcommand{\kphi}{{K\langle \Phi \rangle}}
\newcommand{\kz}{{K\langle z \rangle}}
\newcommand{\kzz}{{K\langle \langle z \rangle\rangle}}
\newcommand{\ktz}{{K[t]\langle z \rangle}}

\newcommand{\kttzz}{{K[[t]]\langle \langle z \rangle\rangle}}

\newcommand{\BQ}{\begin{eqnarray}}
\newcommand{\EQ}{\end{eqnarray}}
\newcommand{\BQn}{\begin{eqnarray*}}
\newcommand{\EQn}{\end{eqnarray*}}

\newcommand{\lb}{\left[}
\newcommand{\rb}{\right]}
\newcommand{\lp}{\left(}
\newcommand{\rp}{\right)}

\newcommand{\fr}{\frac}
\newcommand{\pz}{\frac{\p}{\p z}}

\title[Differential Operator Specializations of NCSF$\mbox{s}$]
{Differential Operator Specializations of Noncommutative Symmetric Functions}

    \author{Wenhua Zhao}      
    \date{\today}
    \begin{document}

\begin{abstract} 
Let $K$ be any unital commutative $\bQ$-algebra and 
$z=(z_1, ... , z_n)$ commutative 
or noncommutative free variables. 
Let $t$ be a formal parameter which commutes 
with $z$ and elements of $K$. We denote 
uniformly by $\kzz$ and 
$\kttzz$ the formal power 
series algebras of $z$ 
over $K$ and $K[[t]]$, respectively.
For any $\alpha \geq 1$, 
let $\cDazz$ be the unital algebra generated by 
the differential operators of $\kzz$ 
which increase the degree in $z$ by at least $\alpha-1$ 
and $ \ataz $ the group of automorphisms $F_t(z)=z-H_t(z)$ 
of $\kttzz$ with $o(H_t(z))\geq \alpha$ and $H_{t=0}(z)=0$.
First, for any fixed $\alpha \geq 1$ and $F_t\in \ataz$, 
we introduce five sequences of differential operators 
of $\kzz$ and show that their generating functions form 
a \cNcs (noncommutative symmetric) system (\cite{GTS-I}) 
over the differential algebra $\cDazz$. 
Consequently, by the universal property of 
the \cNcs system formed by the generating functions of 
certain NCSFs (noncommutative symmetric functions)
first introduced in \cite{G-T},  
we obtain a family of Hopf algebra homomorphisms
$\cS_{F_t}: \cNsf \to \cDazz$ $(F_t\in \ataz)$, 
which are also grading-preserving  when 
$F_t$ satisfies certain conditions.
Note that, the homomorphisms 
$\cS_{F_t}$
above can also be viewed as specializations 
of NCSFs by the differential operators 
of $\kzz$. Secondly, we show that, 
in both commutative and noncommutative cases,
this family $\cS_{F_t}$ 
(with all $n\geq 1$ and $F_t\in \ataz$) 
of differential operator specializations 
can distinguish any two different NCSFs.
Some connections of the results above with 
the quasi-symmetric functions 
(\cite{Ge}, \cite{MR}, \cite{St2})
are also discussed.
\end{abstract}

   

 \bibliographystyle{alpha}
    \maketitle


\renewcommand{\theequation}{\thesection.\arabic{equation}}
\renewcommand{\therema}{\thesection.\arabic{rema}}
\setcounter{equation}{0}
\setcounter{rema}{0}
\setcounter{section}{0}

\section{\bf Introduction}\label{S1}

Let $K$ be any  
unital commutative $\bQ$-algebra and $A$ a 
unital associative but not necessarily 
commutative $K$-algebra. 
Let $t$ be a formal central parameter, 
i.e. it commutes with all elements of $A$, 
and $A[[t]]$ the $K$-algebra 
of formal power series 
in $t$ with coefficients in $A$.
A {\it \cNcs system} over $A$ 
(see Definition \ref{Main-Def}) 
by definition is a $5$-tuple 
$\Omega\in A[[t]]^{\times 5}$ 
which satisfies the defining equations 
(see Eqs.\,$(\ref{UE-0})$--$(\ref{UE-4})$) 
of the NCSFs (noncommutative symmetric functions) 
first introduced and studied 
in the seminal paper \cite{G-T}. 
When the base algebra 
$K$ is clear in the context,
the ordered pair $(A, \Omega)$ 
is also called a {\it \cNcs system}.    
In some sense, a \cNcs system over 
an associative $K$-algebra can be viewed 
as a system of analogs in $A$
of the NCSFs defined by 
Eqs.\,$(\ref{UE-0})$--$(\ref{UE-4})$.
For some general discussions on 
the \cNcs systems, see \cite{GTS-I}.
For a \cNcs  system 
over the Grossman-Larson 
Hopf algebra (\cite{GL}, \cite{F})
of labeled rooted trees, 
see \cite{GTS-IV}.
For more studies on NCSFs, 
see \cite{T}, \cite{NCSF-II}, 
\cite{NCSF-III}, \cite{NCSF-IV}, 
\cite{NCSF-V} and \cite{NCSF-VI}.

One immediate but probably the most 
important example of the \cNcs systems 
is $(\cNsf, \Pi)$ formed 
by the generating functions of 
the NCSFs defined in \cite{G-T}
by Eqs.\,$(\ref{UE-0})$--$(\ref{UE-4})$ 
over the free $K$-algebra $\cNsf$ of NCSFs  
(see Section \ref{S2}). 
It serves as the universal \cNcs system 
over all associative $K$-algebra 
(see Theorem \ref{Universal}).  
More precisely, for any \cNcs system $(A, \Omega)$, 
there exists a unique $K$-algebra homomorphism 
$\cS: \cNsf \to A$ such that
$\cS: \cNsf \to A$ such that 
$\cS^{\times 5}(\Pi) = \Omega$ 
(here we have extended the homomorphism
$\cS$ to $\cS: \cNsf[[t]] \to A[[t]]$ 
by the base extension). 

The universal property of 
the \cNcs system $(\cNsf, \Pi)$
can be applied as follows 
when a \cNcs system $(A, \Omega)$ 
is given.
Note that, as an important topic in 
the symmetric function theory, 
the relations or polynomial identities 
among various NCSFs have been 
worked out explicitly 
(see \cite{G-T}). 
When we apply the $K$-algebra 
homomorphism $\cS:\cNsf \to A$ 
guaranteed by the universal property of 
the system $(\cNsf, \Pi)$ 
to these identities, 
they are transformed into 
identities among the corresponding
elements of $A$ in the system $\Omega$.
This will be a very effective way 
to obtain identities for 
certain elements of $A$ if 
we can show they are 
involved 
in a \cNcs system over $A$.
On the other hand, 
if a \cNcs system 
$(A, \Omega)$ has already been 
well-understood, 
the $K$-algebra homomorphism 
$\cS:\cNsf  \to A$ in turn 
provides a 
a {\it specialization} 
or {\it realization} (\cite{G-T}, \cite{St2})
of NCSFs, 
which may provide some new understandings on NCSFs.
For more studies onthe specializations 
of NCSFs, see the references quoted 
above for NCSFs.

In this paper, motivated by the studies 
on the deformations of formal 
analytic maps of affine spaces 
in \cite{BurgersEq} and \cite{NC-IVP},
we first construct 
a family of \cNcs systems 
over differential operator algebras 
and then study some properties of 
the resulting  specializations 
of NCSFs. 
To be more precise, let us first fix 
the following notations.
Let $z=(z_1, z_2, ... ,z_n)$
be commutative or noncommutative 
free variables and $t$ 
a central parameter, i.e. 
$t$ commutes with $z$.
To keep notation simple, 
we use the notations 
for noncommutative variables 
uniformly for both commutative 
and noncommutative variables $z$. 
Let $\kzz$ (resp.\,\,$\kz$) 
the algebra of formal power series 
(polynomials) in $z$ over $K$. 
For any $\alpha \geq 1$, 
let $\cDazz$ (resp.\,\,$\cDaz$)
be the unital algebra generated by 
the differential operators of $\kzz$ 
(resp.\,\,$\kz$) which increase the degree in $z$ 
by at least $\alpha-1$ 
and $ \ataz $ the group of automorphisms $F_t(z)=z-H_t(z)$ 
of $\kttzz$ with $o(H_t(z))\geq \alpha$ and $H_{t=0}(z)=0$.

First, for any fixed $F_t\in \ataz$, 
we consider the following differential operators,
which are the differential operators
involved in the Taylor series expansions of 
$u(F_t(z))$ and $u(F^{-1}_t(z))$ ($u(z)\in \kzz$), 
the differential operators directly 
related with the D-Log 
(see \cite{E1}--\cite{E3}, \cite{N}, \cite{Z-exp}, 
and \cite{WZ} for the commutative case) 
of $F_t$ and, finally, two sequences of 
differential operators that appeared 
in \cite{NC-IVP} in the study 
of deformations of 
the automorphisms of $\kzz$. 
We show that the generating 
functions of these five sequences 
of differential operators 
form a \cNcs system 
$\Omega_{F_t}$ over $\cDazz$ 
(see Theorem \ref{S-Correspondence}).
Consequently, by the universal properties 
of the \cNcs system
$(\cNsf, \Pi)$, 
we obtain a {\it differential operator 
specialization} $\cS_{F_t}: \cNsf \to \cDazz$, 
which can be shown is also a homomorphism 
of $K$-Hopf algebras.

Secondly, we prove the following
properties of the differential 
operator specializations
$\cSft: \cNsf\to \cDazz$ 
$(F_t\in \ataz)$ above.
We first show in Proposition \ref{cS-graded}
that, for any $F_t\in \ataz$, 
the specialization $\cSft$ is
a homomorphism of graded $K$-Hopf algebras from $\cNsf$ 
to the subalgebra $\cDaz\subset\cDazz$
if and only if 
$F_t(z)=t^{-1}F(tz)$ for some automorphism $F(z)$ of $\kzz$.
Consequently, for any $F_t$ satisfying 
the condition above,  by taking the graded duals, 
we get a graded $K$-Hopf algebra homomorphism
$\cSft^*: \cDaz^*  \to \cQf$ from the
graded dual $\cDaz^*$ 
of $\cDaz$ to the Hopf algebra $\cQf$ 
(\cite{Ge}, \cite{MR} and \cite{St2})
of quasi-symmetric functions 
(see Corollary \ref{S-cQsf}). 
We then show in Theorem \ref{Isomorphism} that, 
with a properly defined group product 
for the set $\text{\bf Hopf}_K(\cNsf, \cDazz)$
of all $K$-Hopf   algebra homomorphisms from 
$\cNsf$, the correspondence $F_t\in \ataz$ to 
$\cSft\in \text{\bf Hopf}_K(\cNsf, \cDazz)$
gives an isomorphism of groups.
Finally, in Theorem \ref{StabInjc}, 
we show that, the family of 
the specializations $\cS_{F_t}$ 
with all $n\geq 1$ 
(note that $n$ is the number of free variables $z_i$) 
and all $F_t\in \ataz$ 
can distinguish any two NCSFs.
 
The arrangement of the paper is as follows.
We first in Section 
\ref{S2} recall the definitions of 
the \cNcs systems and the universal 
\cNcs system $(\cNsf, \Pi)$ 
from NCSFs. In Subsection \ref{S3.1}, 
we first fixed some notations 
and prove some lemmas on the differential 
operators in commutative or noncommutative variables. 
In Subsection \ref{S3.2}, for each fixed $\alpha\geq 1$ 
and $F_t\in \ataz$, we introduce five sequences of
differential operators and show that their generating 
functions actually form a \cNcs system $\Oft$ 
over the differential operator algebra $\cDazz$.
Consequently, we get a differential operator 
$\cS_{F_t}:\cNsf \to \cDazz$, which is also 
a $K$-Hopf algebra. In Subsection \ref{S4}, 
we prove the properties of the differential 
operator specialization 
$\cS_{F_t}:\cNsf \to \cDazz$
that have been explained in the previous 
paragraph.

Finally, some remarks are as follows.  
This paper is the second  of a sequence of papers 
on \cNcs systems over differential operator algebras 
in commutative or noncommutative variables 
and the Grossman-Larson Hopf algebra of
labeled rooted trees as well as their applications 
to NCSF specializations and the inversion problem.
In the followed paper \cite{GTS-IV}, for any nonempty 
$W\subseteq \bN^+$,
a \cNcs system $\Omega_\bT^W$ over 
the Grossman-Larson Hopf algebra 
$\cH^W_{GL}$ (\cite{GL}, \cite{F}) 
of the $W$-labeled rooted trees will be constructed. 
The relations of the \cNcs system 
$(\cH^W_{GL}, \Omega_\bT^W)$
with the \cNcs systems $(\cDazz, \Omega_{F_t})$ 
constructed in this paper 
will be studied in the followed paper \cite{GTS-V}.
In particular, Theorem \ref{StabInjc} derived 
in this paper will be improved to the much smaller 
family of specializations $\cSft: \cNsf\to \cDaz$ with
all $n\geq 1$ and $F_t=z-H_t(z) \in \ataz$
such that, $H_t(z)$ is homogeneous and 
the Jacobian matrix $JH_t$ is strictly lower triangular.
But the proof there is based on Theorem \ref{StabInjc} 
itself and some connections derived in \cite{GTS-V}
among the \cNcs system $(\cNsf, \Pi)$, 
$(\cH^W_{GL}, \Omega_\bT^W)$ 
and $(\cDazz, \Omega_{F_t})$. 
Finally, by the gadget mentioned in the second 
paragraph of this Introduction, by applying 
the specializations $\cSft:\cNsf\to \cDazz$ to the identities 
of the NCSFs in the \cNcs system $(\cNsf, \Pi)$, 
we obtain a host of identities for the differential operators 
in the \cNcs system $(\cDazz, \Oft)$. 
Some of these identities and their consequences 
to the inversion problem (\cite{BCW}, \cite{E})
will be studied in the followed paper 
\cite{GTS-III}. Some other consequences of 
the \cNcs systems $(\cH^W_{GL}, \Omega_\bT^W)$ 
and $(\cDazz, \Omega_{F_t})$  
to the inversion problem will 
also be derived in \cite{GTS-V}.

\renewcommand{\theequation}{\thesection.\arabic{equation}}
\renewcommand{\therema}{\thesection.\arabic{rema}}
\setcounter{equation}{0}
\setcounter{rema}{0}

\section{\bf The Universal \cNcs System from 
Noncommutative Symmetric Functions} \label{S2}

In this section, we first 
recall the definition 
of the \cNcs systems (\cite{GTS-I})
over associative algebras 
and some of the NCSFs 
(noncommutative symmetric functions) 
first introduced and studied 
in the seminal paper \cite{G-T}. 
We then discuss the universal property 
of the \cNcs system formed by the generating 
functions of these NCSFs. The main result that 
we will need later is Theorem \ref{Universal} which 
was proved in \cite{GTS-I}.

Let $K$ be any unital commutative $\bQ$-algebra and 
$A$ any unital associative but not necessarily commutative 
$K$-algebra. Let $t$ be a formal central parameter, 
i.e. it commutes with all elements of $A$, and $A[[t]]$ 
the $K$-algebra of formal power series 
in $t$ with coefficients in $A$. 
First let us recall the following notion formulated 
in \cite{GTS-I}.

\begin{defi} \label{Main-Def}
For any unital associative $K$-algebra $A$, a $5$-tuple $ \Omega=$ 
$( f(t)$, $g(t)$, $d\,(t)$, $h(t)$, $m(t) ) 
\in A[[t]]^{\times 5}$ is said 
to be a {\it \cNcs $($noncommutative symmetric$)$ system}
over $A$ if the following equations are satisfied.
\begin{align}
&f(0)=1 \label{UE-0}\\
& f(-t)  g(t)=g(t)f (-t)=1, \label{UE-1}   \\
& e^{d\,(t)} = g(t), \label{UE-2} \\
& \frac {d g(t)} {d t}= g(t) h(t), \label{UE-3}\\ 
& \frac {d g(t)}{d t} =  m(t) g(t).\label{UE-4}
\end{align}
\end{defi}

When the base algebra $K$ is clear in the context, we also call 
the ordered pair $(A, \Omega)$ a {\it \cNcs system}. 
Since \cNcs systems often come from generating functions 
of certain elements of $A$ that are under the consideration, 
the components of $\Omega$ will also be refereed as 
the {\it generating functions} of their coefficients. 

All $K$-algebras $A$ that we are going to work on in this paper
are $K$-Hopf algebras. We will freely use some standard results 
from the theory of bi-algebras and Hopf algebras, 
whose proofs can be found in the standard text books 
\cite{Abe}, \cite{Knu} and \cite{Mon}.

The following result proved in \cite{GTS-I} later
will be useful to us.

\begin{propo}\label{bialg-case}
Let $(A, \Omega)$ be a \cNcs system as above. 
Suppose $A$ is further a $K$-bialgebra. 
Then the following statements are equivalent. 
\begin{enumerate}
\item[$(a)$] The coefficients of $f(t)$ form a sequence of divided powers of $A$.

\item[$(b)$] The coefficients of $g(t)$ form a sequence of divided powers of $A$.

\item[$(c)$] One $(\text{hence also all})$ 
of \,$d(t)$, $h(t)$ and $m(t)$ has all 
its coefficients primitive in $A$.
\end{enumerate}
\end{propo}

Next, let us recall some of the NCSFs 
first introduced and studied in (\cite{G-T}). 

Let $\Lambda=\{ \Lambda_m\,|\, m\geq 1\}$ 
be a sequence of noncommutative 
free variables and $\cNsf$ or $\klam$ 
the free associative algebra 
generated by $\Lambda$ over $K$.  
For convenience, we also set $\Lambda_0=1$.
We denote by
$\lambda (t)$ the generating function of 
$\Lambda_m$ $(m\geq 0)$, i.e. we set
\begin{align}
\lambda (t):= \sum_{m\geq 0} t^m \Lambda_m 
=1+\sum_{k\geq 1} t^m \Lambda_m.
\end{align}

In the theory of NCSFs (\cite{G-T}), 
$\Lambda_m$ $(m\geq 0)$ is 
the noncommutative analog
of the $m^{th}$ classical (commutative) 
elementary symmetric function 
and is called the {\it $m^{th}$ 
$(\text{noncommutative})$ 
elementary symmetric function.}

To define some other NCSFs, we consider 
Eqs.\,$(\ref{UE-1})$--$(\ref{UE-4})$ 
over the free $K$-algebra $\cNsf$
with $f(t)=\lambda(t)$. The 
solutions for $g(t)$, $d\,(t)$, 
$h(t)$, $m(t)$ exist and are unique, 
whose coefficients will be the NCSFs 
that we are going to define.
Following the notation in \cite{G-T} 
and \cite{GTS-I}, we denote the resulting  
$5$-tuple by 
\begin{align}
\Pi:= (\lambda(t),\, \sigma(t),\, \Phi(t),\, \psi(t),\, \xi(t))
\end{align}
and write the last 
four generating functions of 
$\Pi$ explicitly as follows.

\allowdisplaybreaks{
\begin{align}
\sigma (t)&=\sum_{m\geq 0} t^m S_m,  \label{lambda(t)} \\
\Phi (t)&=\sum_{m\geq 1} t^m \frac{\Phi_m}m  \label{Phi(t)}\\
\psi (t)&=\sum_{m\geq 1} t^{m-1} \Psi_m, \label{psi(t)}\\
\xi (t)&=\sum_{m\geq 1} t^{m-1} \Xi_m.\label{xi(t)}
\end{align}}

Note that, by Definition \ref{Main-Def},
the $5$-tuple $\Pi$ defined above is just 
the unique \cNcs system 
with $f(t)=\lambda(t)$ in 
Eq.\,(\ref{lambda(t)}) 
over the free $K$-algebra 
$\cNsf$.

Following \cite{G-T},
we call $S_m$ ($m\geq 1$) the 
{\it $m^{th}$ $(\text{noncommutative})$ complete 
homogeneous symmetric function} and
$\Phi_m $ (resp.\,\,$\Psi_m$) 
the {\it $m^{th}$ power sum symmetric function 
of the second $($resp.\,\,first$)$ kind}. 
Following \cite{GTS-I}, 
we call $\Xi_m \in \cNsf$ $(m\geq 1)$ 
the {\it $m^{th}$ $(\text{noncommutative})$ 
power sum symmetric function of the third kind}.

The following two propositions proved in \cite{G-T} 
and \cite{NCSF-II} will be very useful 
for our later arguments.

\begin{propo}\label{bases}
For any unital commutative $\bQ$-algebra $K$,
the free algebra $\cNsf$ is freely generated
by any one of the families of the NCSFs 
defined above.
\end{propo}

\begin{propo}\label{omega-Lambda}
Let $\omega_\Lambda$ be the anti-involution of 
$\cNsf$ 
which fixes $\Lambda_m$ $(m\geq 1)$.
Then, for any $m\geq 1$, we have
\begin{align}
\omega_\Lambda (S_m)&=S_m,  \label{omega-Lambda-e1}\\
\omega_\Lambda (\Phi_m)&=\Phi_m,\label{omega-Lambda-e2} \\
\omega_\Lambda (\Psi_m)&=\Xi_m. \label{omega-Lambda-e3}
\end{align}
\end{propo}

Next, let us recall the following graded 
$K$-Hopf algebra structure 
of $\cNsf$. It has been shown in 
\cite{G-T} that $\cNsf$ is the universal enveloping algebra 
of the free Lie algebra generated 
by $\Psi_m$ $(m\geq 1)$. Hence, it has a $K$-Hopf  
algebra structure as all other universal enveloping algebras 
of Lie algebras do. Its co-unit $\epsilon:\cNsf \to K$,
 co-product $\Delta$ and 
 antipode $S$ are uniquely determined by 
\begin{align}
\epsilon (\Psi_m)&=0, \label{counit} \\
\Delta (\Psi_m) &=1\otimes \Psi_m +\Psi_m\otimes 1, \label{coprod}\\
S(\Psi_m) & =-\Psi_m,\label{antipode}
\end{align}
for any $m\geq 1$. 

Next, we introduce the {\it weight} of NCSFs 
by setting the weight of 
any monomial $\Lambda_{m_1}^{i_1} 
\Lambda_{m_2}^{i_2} \cdots \Lambda_{m_k}^{i_k}$
to be $\sum_{j=1}^k i_j m_j$. 
For any $m\geq 0$, we denote by $\cNsf_{[m]}$ 
the vector subspace of $\cNsf$ spanned 
by the monomials of $\Lambda$ 
of weight $m$. Then it is easy to see that 
\begin{align}\label{Grading-cNsf}
\cNsf=\bigoplus_{m\geq 0} \cNsf_{[m]}, 
\end{align}
which provides a grading for $\cNsf$. 

Note that, it has been shown in \cite{G-T}, 
for any $m\geq 1$, the NCSFs 
$S_m, \Phi_m, \Psi_m \in  \cNsf_{[m]}$. 
By Proposition \ref{omega-Lambda}, 
this is also true for the NCSFs $\Xi_m$'s.
By the facts above and 
Eqs.\,(\ref{counit})--(\ref{antipode}), 
it is also easy to check that, 
with the grading given in Eq.\,(\ref{Grading-cNsf}), 
$\cNsf$ forms a graded $K$-Hopf algebra. 
Its graded dual is given 
by the space $\cQf$ of quasi-symmetric functions, 
which were first introduced by I. Gessel \cite{Ge} 
(see \cite{MR} and \cite{St2} for more discussions).

Now we come back to our discussions on the \cNcs systems. 
From the definitions of the NCSFs above, 
we see that $(\cNsf, \Pi)$ obviously forms a \cNcs system.
More importantly, as shown in Theorem $2.1$ in \cite{GTS-I}, 
we have the following important theorem on 
the \cNcs system $(\cNsf, \Pi)$. 

\begin{theo}\label{Universal}
Let $A$ be a $K$-algebra and $\Omega$ 
a \cNcs system over $A$. Then, 

$(a)$ There exists a unique $K$-algebra homomorphism 
$\cS: \cNsf\to A$ such that 
$\cS^{\times 5} (\Pi)=\Omega$.

$(b)$ If $A$ is further  a $K$-bialgebra $($resp.\,\,$K$-Hopf algebra$)$ 
and one of the equivalent statements in Proposition \ref{bialg-case} 
holds for the \cNcs system $\Omega$, then $\cS: \cNsf\to A$ is also 
a homomorphism of $K$-bialgebras $($resp.\,\,$K$-Hopf algebras$)$.
\end{theo}

\begin{rmk}\label{Comm-Case}
By applying the similar arguments as in the 
proof of Theorem \ref{Universal}, 
or simply taking the quotient over 
the two-sided ideal generated by the commutators 
of $\Lambda_m$'s, it is easy to see that, 
over the category of commutative $K$-algebras, 
the universal \cNcs system 
is given by the generating functions of
the corresponding classical 
$($commutative$)$ symmetric functions $($\cite{Mc}$)$.
\end{rmk}

\begin{rmk}\label{Motivation}
One direct consequence of Theorem \ref{Universal} 
above is as follows. 
Note that the relations or polynomial identities 
between any two families 
of NCSFs in the first four components of 
$\Pi$ have been given explicitly in \cite{G-T}. 
By applying the anti-automorphism $\omega_\Lambda$
in Proposition $\ref{omega-Lambda}$, 
one can easily derive the relations of 
the NCSFs $\Xi_m$'s with 
other NCSFs in $\Pi$ $($for example, 
see $\S 4.1$ in \cite{GTS-III} 
for a complete list$)$.
By applying the homomorphism $\cS$, 
we get a host of identities among 
the corresponding elements of $A$. 
This will be a very effective 
method to prove identities for 
the elements of $A$ which are 
involved in a \cNcs system 
over $A$. On the other hand, 
if the \cNcs system 
$(A, \Omega)$ has already
been well understood, the homomorphism 
$\cS: \cNsf\to A$, usually called 
a {\it specialization} of NCSFs, 
can also be used to study certain 
properties of NCSFs.
\end{rmk}

\renewcommand{\theequation}{\thesection.\arabic{equation}}
\renewcommand{\therema}{\thesection.\arabic{rema}}
\setcounter{equation}{0}
\setcounter{rema}{0}

\section{\bf \cNcs Systems over Differential Operator Algebras}\label{S3}

Let $K$ be any unital commutative $\bQ$-algebra
and $z=(z_1, z_2, ... , z_n)$ free variables, i.e.
commutative or noncommutative independent 
variables.
Let $t$ be a formal central parameter, 
i.e. it commutes with $z$ 
and all elements of $K$.
We denote by $\kzz$ 
 and 
$\kttzz$ the $K$-algebras of 
formal power series in $z$ over  
$K$ and $K[[t]]$, 
respectively.\footnote{Since most of the results 
as well as their proofs in this paper
do not depend on the commutativity 
of the free variables $z$,  we will 
not distinguish the commutative 
and the noncommutative case, 
unless stated otherwise, 
and adapt the notations 
for noncommutative variables 
uniformly for the both cases.}
For any positive integer $\alpha \geq 1$, 
we let $\cDazz$ denote the $K$-algebra 
of the differential operators of $\kzz$ 
which increase the degree in $z$ 
by at least $\alpha-1$. We also 
fix the notation $\ataz$ for the set
of all the automorphisms $F_t(z)$ of 
$\kttzz$ over $K[[t]]$ 
which have the form
$F(z)=z-H_t(z)$ for some 
$H_t(z)\in \kttzz^{\times n}$ 
with $o(H_t(z))\geq \alpha$ 
and $H_{t=0}(z)=0$.
It is easy to check that 
$\ataz$ actually forms 
a subgroup of 
the automorphism group 
of $\kttzz$.

In Subsection \ref{S3.1}, we  
fix more notation and 
prove some simple results 
on the differential operators 
in commutative or noncommutative 
free variables $z$.
In Subsection \ref{S3.2}, 
for any fixed automorphism 
$F_t(z)\in \ataz$,  
we introduce five families of differential 
operators associated with  $F_t(z)$ 
and its inverse $G_t(z):=F_t^{-1}(z)$
and show that their generating 
functions form a \cNcs system $\Omega_{F_t}$ 
over the differential operator algebra $\cDazz$. 
Consequently, by the universal
property (see Theorem \ref{Universal})
of the \cNcs system 
$(\cNsf, \Pi)$ from NCSFs, 
we obtain a family of differential operator 
specializations $\cS_{F_t}: \cNsf\to \cDazz$ 
$(F_t(z)\in \ataz)$ for NCSFs. 
More properties of the specializations 
$\cS_{F_t}$ $(F_t\in \ataz)$ will 
be studied in next section.

\subsection{Differential Operators in 
(Noncommutative) Free Variables}\label{S3.1}

Let $K$, $z$ and $t$ as fixed above. 
In this subsection, 
we mainly fix 
more notations and prove some simple results 
on the differential operators in $z$.
All the results proved in this section 
should be well-known, especially 
in the commutative case. 
But for the completeness, 
especially when the noncommutative 
case is concerned, 
we also include proofs here. 
As we mentioned early, 
in this subsection 
as well as in the rest of this paper,
we do not assume the commutativity 
of our free variables 
$z$, unless stated otherwise.
 
Recall that, a {\it $K$-derivation} or simply a {\it derivation} 
of $\kzz$ is a $K$-linear map $\delta: \kzz\to \kzz$ 
which satisfies the Leibnitz rule,
i.e. for any $f, g\in \kzz$, we have
 \begin{align}\label{Leibnitz}
\delta (fg)=(\delta f)g+f(\delta g).
\end{align}
We denote by 
$\cDrkzz$ or $\cDrzz$, when the base algebra
$K$ is clear in the context, 
the set of all $K$-derivations 
of $\kzz$.
The unital subalgebra of 
$\text{End}_K (\kzz)$ 
(the set of endomorphisms of $\kzz$ as a $K$-vector space, 
not as a $K$-algebra) 
generated by all
$K$-derivations 
of $\kzz$ is denoted by 
$\cDkzz$ or $\cDzz$. 
Elements of 
$\cDzz$ are called the
{\it $(\text{formal})$ 
differential operators} 
in the free variables $z$.

For any $\alpha \geq 1$, we denote by 
$\cDrazz$ the set of 
$K$-derivations of $\kzz$ 
which increase 
the degree in $z$ by 
at least $\alpha-1$. 
The unital subalgebra of 
$\cDzz$ generated by elements of 
$\cDrazz$ will be denoted by $\cDazz$. 
Note that, by the definitions above,
the operators of scalar multiplications
are in $\cDzz$ and $\cDazz$, 
but not in $\cDrazz$.
When the base algebra is $K[[t]]$ 
instead of $K$ itself,
the corresponding notation 
$\cDrzz$, $\cDzz$, $\cDrazz$ 
and $\cDazz$ will be denoted by  
$\cDrtzz$, $\cDtzz$, $\cDrtazz$ 
and $\cDtazz$, respectively.
For example,  $\cDrtazz$ stands for 
the set of all $K[[t]]$-derivations of the $K[[t]]$-algebra
$\kttzz$ which increase the degree in $z$ 
by at least $\alpha-1$. 
Note that, $\cDrtazz=\cDrazz[[t]]$ and 
$\cDtazz=\cDazz[[t]]$.

For any $1\leq i\leq n$ and $u(z)\in \kzz$, 
we denote by $\lb u(z) \fr \p{\p z_i}\rb $
\footnote{The reason we put a bracket $[\cdot]$ 
in the notation for this derivation of $\kzz$ is to avoid 
any possible confusion caused by a subtle point described 
in the {\bf Warning} below.}
the $K$-derivation which maps $z_i$ to $u(z)$ and $z_j$ to $0$ 
for any $j\neq i$.
For any $\Vec{u}=(u_1, u_2, \dots , u_n)\in \kzz^{\times n}$, 
we set 
\begin{align}\label{Upz}
[\Vec{u}\pz]:=\sum_{i=1}^n [u_i \fr\p{\p z_i}].  
\end{align}

{\bf Warning:} 
{\it When $z$ are noncommutative free variables, 
we in general do {\bf not} have
$\lb u(z) \fr \p{\p z_i}\rb  g(z)  = u(z)  
\fr {\p g}{\p z_i}$ for all 
$u(z), g(z)\in \kzz$.
For example, let $g = z_j z_i$ with $j\neq i$, we have
\begin{align*}
[u \fr \p{\p z_i}](z_jz_i)& = ([u \fr \p{\p z_i}]z_j )z_i+z_j 
([u \fr \p{\p z_i}]z_i)= z_j u(z),\\
u(z) \fr {\p g}{\p z_i} & =u(z)z_j,
\end{align*}
which are not equal to each other unless $u(z)$ commutes with $z_j$.}

\vskip3mm

With the notation above, 
it is easy to see that any 
$K$-derivation $\delta$ 
of $\kzz$ can be written uniquely  as 
$\sum_{i=1}^n \lb f_i(z)\fr\p{\p z_i}\rb$ 
with $f_i(z)=\delta \cdot z_i\in \kzz$ 
$(1\leq i\leq n)$.
Also, as in the commutative case, 
we have the following Lie bracket 
relation in the noncommutative case,
namely, for any derivations 
$\delta=\lb \vec{u} \pz\rb$ and 
$\eta=\lb \vec{v} \pz \rb$ with 
$\vec{u}, \vec{v} \in \kzz^{\times n}$, 
we have
\begin{align}
[\delta, \eta]=\lb (\delta \vec{v}) \pz \rb -\lb (\eta \vec{u})\pz \rb,
\end{align}
where $[\delta, \eta]$ is the commutator of $\delta$ and $\eta$.

With the bracket above, 
$\cDrzz$ forms a Lie algebra 
and its universal enveloping algebra is exactly 
the differential operator algebra
$\cDzz$.  Consequently, $\cDzz$ has a  
$K$-Hopf algebra structure as all 
other enveloping algebras of Lie algebras do.
In particular,
Its coproduct $\Delta$, antipode $S$ and counit $\epsilon$
are respectively determined by the following properties:
for any $\delta\in \cDrzz$, 
\begin{align}
\Delta(\delta)&= 1\otimes\delta+\delta\otimes 1,\label{Coprd-delta} \\
S(\delta)&=-\delta, \label{antipd-delta}\\
\epsilon (\delta) &=\delta \cdot 1. \label{Counit-delta}
\end{align}

Next, we introduce the following two 
operations for the $K$-derivations of $\kzz$.

First, for any $\phi, \, \delta \in \cDrzz$ with
$\delta=\lb \vec{u}\pz\rb$ for some $\vec{u} \in \kzz^{\times n}$, 
we set
\begin{align}\label{collapse}
\phi \triangleright \delta:=\lb (\phi \vec{u})(z)\pz\rb. 
\end{align}

Secondly, for any $K$-derivations 
$\delta_i=\lb \Vec v_i(z)\pz \rb$ $(1\leq i\leq m)$
with $\Vec v_i(z)\in \kzz^{\times n}$, we define 
a new linear operator 
$B_+(\delta_1, \delta_2, \dots , \delta_m)$ as follows. 
Let $w=(w_1, w_2, \dots , w_n)$ be 
another $n$ free variables which are 
independent and do not commute 
with the free variables $z$. 
We define $B_+(\delta_1, \delta_2, \dots , \delta_m)$ 
by setting, for any $u(z)\in \kzz$, 
\begin{align}\label{Def-B+}
& B_+(\delta_1, \delta_2, \dots , \delta_m)u(z):=\\
&\quad\quad\quad\quad
\left. \lb \Vec v_1(w)\pz\rb  \lb \Vec v_2(w)\pz \rb\cdots 
\lb \Vec v_m(w)\pz \rb u(z)\, \right |_{w=z}.\nno 
\end{align}

Note that $B_+(\delta_1, \delta_2, \dots , \delta_m)$ 
is multi-linear and symmetric in the components $\delta_i$ 
$(1\leq i\leq m)$. When $m=1$, we have $B_+(\delta_1)=\delta_1$.

Furthermore, for any $k_i\geq 0$ $(1\leq i\leq m)$, 
we let 
$B_+(\delta_1^{k_1}, \delta_2^{k_2}, \dots , \delta_m^{k_m})$ 
denote the operator obtained by applying $B_+$ 
to the multi-set of $j_1$-copies of $\delta_1$; 
$j_2$-copies of $\delta_2$, ..., 
$j_m$-copies of $\delta_m$.

Next, we show that 
$B_+(\delta_1, \delta_2, \dots , \delta_m)$ 
is still a differential operator in $z$, i.e. 
$B_+(\delta_1, \delta_2, \dots , \delta_m) \in \cDzz $. 
But first we need prove the following lemma.

\begin{lemma}\label{L2.1.2}
For any  $\phi, \, \delta_i \in \cDrzz$ $(1\leq i\leq m)$, we have
\begin{align}\label{L2.1.2-e1}
\phi \cdot B_+(\delta_1, \delta_2, \dots , \delta_m)&=
B_+(\phi, \delta_1, \delta_2, \dots , \delta_m) \\* 
& \quad \quad \quad \quad  +\sum_{i=1}^m B_+
(\delta_1, \dots , \phi\triangleright \delta_i, \dots   \delta_m).\nno
\end{align}
\end{lemma}

\pf First, since
$B_+(\delta_1, \delta_2, \dots , \delta_m)$ 
is multi-linear in the components
$\delta_i$'s, we may assume  
$\delta_i=\lb a_i(z)\frac \p{\p z_{k_i}}\rb $ $(1\leq i\leq m)$ 
for some $a_i(z)\in \kzz$ and $1\leq k_i\leq n$.  
Furthermore, to show Eq.\,$(\ref{L2.1.2-e1})$,  
we only need show its two sides 
have same values at all monomials of $z$.

Now, set $\Psi:=B_+(\delta_1, \delta_2, \dots , \delta_m)$ 
and let $u(z)$ be any monomial of $z$. By
Eq.\,$(\ref{Def-B+})$,  
we know that $\Psi u(z)$ 
is the sum of all the terms 
obtained by replacing $m$-copies 
$z_{k_i}$'s in the monomial $u(z)$
by the corresponding $a_i(z)$'s 
in all possible ways. 
Consequently, each of these terms 
is a monomial in $z_i$'s 
and $m$-copies $a_i(z)$'s. 
Now, we apply the derivation $\phi$  
to $\Psi u(z)$. 
By the Leibnitz rule, 
we know that $\phi$ either 
lands on a variable $z_i$ or 
a copy of $a_{i}(z)$'s. Then,
it is easy to check that 
the sum of all the terms 
obtained in the former case 
is same as $B_+(\phi, \delta_1, \delta_2, \dots , \delta_m)u(z)$;
while, the sum of all the terms 
obtained in the later case 
is same as those obtained
by applying $\sum_{i=1}^m B_+
(\delta_1, \dots , \phi\triangleright \delta_i, \dots   \delta_m)$ 
to $u(z)$.
\epfv

\begin{corol}\label{C2.1.3}
For any $\alpha \geq 1$ and $\delta_i \in \cDrazz$ 
$(1\leq i\leq m)$, we have 
$ B_+(\delta_1, \delta_2, \dots , \delta_m) \in \cDazz$.
\end{corol}
\pf We use the mathematical induction on 
$m\geq 1$. When $m=1$, we have 
$B_+(\delta_1)=\delta_1$, 
hence nothing needs to prove.

Now, let $m\geq 2$. By Eq.\,$(\ref{L2.1.2-e1})$,  
we have
\begin{align*} 
 B_+ (\delta_1, \dots , \delta_m) = 
 \delta_1  B_+(\delta_2, \dots , \delta_m)
-\sum_{i=2}^{m} B_+
(\delta_2, \dots , \delta_1 \triangleright \delta_i, \dots   \delta_{m}).
\end{align*}
Note that, by Eq.\,$(\ref{collapse})$, 
$\delta_1 \triangleright \delta_i$ 
$(2\leq i\leq m)$ is still a derivation 
in $\cDrazz$. Therefore,  
from the equation above 
and the induction assumption, we see that 
$B_+ (\delta_1, \delta_2, \dots , \delta_m)\in \cDazz$.
\epfv

Next, let us consider the Taylor series expansions 
for formal power series in commutative or
noncommutative variables $z$. 
First, we have the following lemma which 
can be proved by a similar argument as 
in the commutative case.

\begin{lemma}\label{L2.1.4}
Let $v=(v_1, v_2, ... , v_n)$ be 
another $n$ free variables 
which are independent and do not commute with $z$. 
Let $t$ be a formal parameter which commutes 
with both $z$ and $v$.
Then, for any $u(z)\in \kzz$, 
we have
\begin{align}\label{L2.1.4-e1}
u(z+tv)= \sum_{k \geq 0} \frac {t^{k}}{k!}
 B_+\left( \lb v \pz \rb^k \right ) \, u(z).
\end{align}
\end{lemma}

\begin{propo}\label{C2.1.5}
For any $\alpha\geq 1$ and $F_t(z)=z-H_t(z)\in \ataz$, 
we have the following Taylor series expansion 
for any $u(z)\in \kzz$.
\begin{align}\label{C2.1.5-e1}
u(F_t(z))= \sum_{k \geq 0} \frac {(-1)^{k}}{k!}
B_+\lp \lb H_t(z) \pz \rb^k \rp \, u(z).
\end{align}
\end{propo}
\pf The proposition follows directly by setting 
$tv=-H_t(z)$ in Eq.\,$(\ref{L2.1.4-e1})$ and then 
using the definition of the operation 
$B_+$ in Eq.\,$(\ref{Def-B+})$. 
\epfv

Finally, let us prove the following lemma which will be needed later.

\begin{lemma}\label{L3.1.5}
Let $z=(z_1, z_2, \dots, z_n)$ be commutative free variables and $\Phi \in \cDzz$.
Suppose that, there exists $N>0$ such that $\Phi \, u(z)=0$ for any $u(z)\in \kzz$ 
with $o(u(z))\geq N$. Then $\Phi=0$. 
\end{lemma}

\pf First, let $\text{End}_K(\kzz)$ be the set of 
$K$-linear maps from $\kzz$ to $\kzz$
and $\mathcal A$ the unital subalgebra of $\text{End}_K(\kzz)$
generated by $K$-derivations and 
the linear operators given by multiplications 
by elements of $\kzz$. Let $\mathcal N$ be the set of all elements 
$\Psi \in \mathcal A$ such that, for some $N\geq 1$ (depending on $\Psi$), 
$\Psi \, u(z)=0$ for any $u(z)\in \kzz$ with $o(u(z))\geq N$. 
It is straightforward to check that $\mathcal N$ forms a left ideal 
of $\mathcal A$. Furthermore, for any polynomial $b(z)\in \kzz$, 
we denote by $L_{b(z)}$ the operator of multiplication by $b(z)$. 
Then, it is easy to check that we also have 
$\mathcal N L_{b(z)}\subset \mathcal N$.

Next we use induction on the {\it order} 
$\text{Ord}\,(\Phi)$ of $\Phi$ to show 
that $\mathcal N=0$. First, if $\text{Ord}\,(\Phi) = 0$, 
then $\Phi$ is just an operator of multiplication 
by an element of $\kzz$, and the lemma obviously holds. 
Assume the lemma holds for all 
$\Phi\in \mathcal N$ with  $\text{Ord}(\Phi) \leq m$. 
Consider the case that $\text{Ord}(\Phi)= m+1$.
Let $\delta:=(\delta_1, \delta_2, \dots , \delta_k)$ 
be a sequence of $K$-derivations such that 
$\Phi$ can be written as
\begin{align}\label{L3.1.5-pe1}
\Phi=\sum_{\substack{I \in \bN^k \\|I|\leq m+1}} a_I(z) \, 
\delta^I 
\end{align}
for some $a_I(z)\in \kzz$. 

Note that, in general, we have
\begin{align}\label{L3.1.5-pe2}
\Phi L_{z_1} 
&=[\Phi, L_{z_1}] + L_{z_1} \Phi ,
\end{align}
where $[\Phi, L_{z_1}]$ denotes the commutator of 
the operators $\Phi$ and $L_{z_1}$.

First, by using the form of $\Phi$ in Eq.\,(\ref{L3.1.5-pe1}),
it is easy to check that 
$\text{Ord}([\Phi, L_{z_1}]) \leq \text{Ord}(\Phi)-1$.
Secondly,  by Eq.\,(\ref{L3.1.5-pe2}) and the facts 
about $\mathcal N$ mentioned in the first paragraph 
of the proof, we have $\Phi L_{z_1}\in \mathcal N$ 
and then $[\Phi, L_{z_1}]\in \mathcal N$.
By our induction assumption, we have
$[\Phi, L_{z_1}]=0$. Therefore, $\Phi$ 
commutes with the left multiplication 
by $z_1$. 
Hence it also commutes with 
the left multiplication by $z_1^m$ for 
any $m\geq 1$. 
Now let $N$ be a positive integer
such that $\Phi\cdot u(z)=0$ for 
any $u(z)\in \kzz$ with $o(u(z))\geq N$.
Then, for any $f(z)\in \kzz$,  we have 
$\Phi( z_1^N f(z) )=z_1^N \Phi f(z)=0$. 
Therefore, $\Phi  f(z)=0$ 
for any $f(z)\in \kzz$. Hence $\Phi=0$.
\epfv

\subsection{\cNcs Systems over Differential Operator Algebras} \label{S3.2}

Let $K$, $z$, $t$, $\cDrazz$ and $\cDazz$ 
$(\alpha\geq 1)$ as fixed 
in the previous subsection. 
Recall that, we also have 
defined $\ataz$ to be the set
of all the automorphisms $F_t(z)$ of 
$\kttzz$ over $K[[t]]$, 
which are of the form
\begin{align}\label{Ft}
F_t(z)=z-H_t(z)
\end{align}
for some $H_t(z)\in \kttzz^{\times n}$ 
with $o(H_t(z))\geq \alpha$ 
and $H_{t=0}(z)=0$. 
Note that, for any $F_t \in \ataz$ as above, 
its inverse map 
can always be written 
uniquely as 
\begin{align}\label{Gt}
G_t(z):=F_t^{-1}(z)=z+M_t(z)
\end{align}
for some $M_t(z)\in \kttzz^{\times n}$ 
with $o(M_t(z))\geq \alpha$ 
and $M_{t=0}(z)=0$. 
Throughout the rest of this section, 
we will fix an arbitrary $F_t\in \ataz$ 
and always let $H_t(z)$, $G_t(z)$ 
and $M_t(z)$ be determined as 
in Eqs.\,$(\ref{Ft})$ 
and $(\ref{Gt})$. 

Note that, $F_t \in \ataz$ can be viewed 
as a deformation parameterized by $t$ of 
the formal map $F(z):=F_{t=1}(z)$, 
when it makes sense. For more studies of 
$F_t\in \ataz$ from the deformation 
point view, see \cite{BurgersEq} 
and \cite{NC-IVP}. Actually,   
the construction of the \cNcs system 
given in this subsection 
is mainly motivated by and also depends 
on the studies of $F_t\in \ataz$ given 
in \cite{BurgersEq} 
and \cite{NC-IVP}.

First, let us introduce the following
five sequences of differential operators 
associated with the fixed automorphism 
$F_t(z)\in \ataz$ and its inverse $G_t(z)$.
As we will see later, the generating functions of 
these differential operators will form a \cNcs system 
over the differential operator algebra $\cDazz$.

The first two sequences of differential operators 
come from $F_t(z)$ and $G_t(z)$ as follows. 

\begin{lemma}\label{TaylorExpansion}
There exist unique sequences $\{\lambda_m | m \geq 0 \}$ 
and $\{s_m | m\geq 0 \}$ of elements of $\cDazz$ 
with $\lambda_0=s_0=1$ such that, for any $u_t(z)\in \kttzz$, we have
\begin{align}
(\sum_{m=0}^\infty (-1)^m t^m \lambda_m)u_t(z)&=u_t(F_t), \label{TaylorExpansion-e1} \\
(\sum_{m=0}^\infty t^m s_m)u_t(z)&=u_t(G_t) \label{TaylorExpansion-e2}.
\end{align}
\end{lemma}

The signs appearing in Eq.\,(\ref{TaylorExpansion-e1}) 
as well as somewhere else 
in this subsection are chosen in such a way that 
the correspondence $\cS_{F_t}$ later 
in Theorem \ref{S-Correspondence} 
between the NCSFs in the universal 
\cNcs system $(\cNsf, \Pi)$
and the differential operators 
defined in this subsection  
will be in the simplest form.

\pf First, let us show the uniqueness. 
By Eqs.\,(\ref{TaylorExpansion-e1}) 
and (\ref{TaylorExpansion-e2}), for any $m\geq 0$ and 
$u(z)\in \kzz \subset \kttzz$,
we have 
\begin{align*}
\lambda_m u(z)& =\res_{t=0}\,\,  u(F_t) \, t^{-m-1} \,\\
s_m u(z)&=\res_{t=0}\, \, u(G_t)\, t^{-m-1} \, 
\end{align*}
Hence, as differential operators of $\kzz$, 
$\lambda_m$ and $s_m$ $(m\geq 1)$ are 
uniquely determined by 
Eqs.\,(\ref{TaylorExpansion-e1}) 
and (\ref{TaylorExpansion-e2}), 
respectively. 

To show the existence, 
we first consider the Eqs.\,$(\ref{TaylorExpansion-e1})$ 
and $(\ref{TaylorExpansion-e2})$ with $u_t(z) \in \kzz$, 
i.e. $u_t(z)$ does not depend on $t$. 
By Corollary \ref{C2.1.5}, for any $u(z)\in \kzz$, 
we have

\begin{align}\label{TaylorExpansion-pe1}
u(F_t(z))= \lp \sum_{k \geq 0} \frac {(-1)^{k}}{k!}
B_+\lp \lb H_t(z) \pz \rb^k \rp \rp \, u(z).
\end{align}

Note that, by Corollary \ref{C2.1.3} and 
the condition $o(H_{t}(z))\geq \alpha$ 
(since $F_t\in \ataz$), 
all the operators involved in the equation above 
are differential operators in 
$\cDtazz$ which is same as 
$\cDazz[[t]]$. 
Therefore, the bracketed sum 
of the differential operators 
in Eq.\,$(\ref{TaylorExpansion-pe1})$ above 
can be written 
as $\sum_{m=0}^\infty (-1)^m t^m \lambda_m$
for some $\lambda_m \in \cDazz$ $(m\geq 0)$. 
With the new differential operators $\lambda_m$'s, 
Eq.\,$(\ref{TaylorExpansion-pe1})$ becomes

\begin{align}\label{TaylorExpansion-pe2}
u(F_t(z))= \lp \sum_{m=0}^\infty (-1)^m t^m \lambda_m \rp u(z).
\end{align}

Note that $H_{t=0}(z)=0$ since $F_t\in \ataz$. 
By setting $t=0$ in Eq.\,(\ref{TaylorExpansion-pe2}) above, 
we get $\lambda_0=1$.

Now we show Eq.\,(\ref{TaylorExpansion-pe2})
also holds for any $u_t(z)\in \kttzz$. 
First, we write 
$u_t(z)=\sum_{k\geq 0} u_k(z)\, t^k$ 
with $u_k(z)\in \kzz$ 
$(k\geq 0)$ and then apply 
Eq.\,(\ref{TaylorExpansion-pe2})
to each $u_k(z)\in \kzz$. Then,
it is easy to see that 
Eq.\,(\ref{TaylorExpansion-pe2}) also holds 
for $u_t(z)\in \kttzz$. 
Therefore, we have proved 
the existence of 
the sequence $\{\lambda_m \in \cDazz\, |\, m\geq 0 \}$.
The existence of the sequence $\{s_m \in \cDazz\, |\, m\geq 0 \}$
can be proved similarly with $F_t$ replaced by $G_t$.
\epfv

We denote by $f(t)$ and $g(t)$  
the generating functions of the differential operators
$\lambda_m$ and $s_m$ $(m\geq 0)$, respectively, i.e. we set 
\begin{align}
f(t):= \sum_{m=0}^\infty t^m \lambda_m, \label{Def-f(t)}\\
g(t) := \sum_{m=0}^\infty t^m s_m. \label{Def-g(t)}
\end{align}

We will also view $f(t)$ and $g(t)$ 
as differential operators
of $\kttzz$ in $\cDazz$. Furthermore, 
for any $u_t(z)\in \kttzz$,  by
Eqs.\,(\ref{TaylorExpansion-e1}) 
and (\ref{TaylorExpansion-e2}), 
we have 
\begin{align}
f(-t)u_t(z)&=u_t(F_t),\label{NewTaylorExpansion-e1}\\
g(t)u_t(z)&=u_t(G_t). \label{NewTaylorExpansion-e2}
\end{align}

From the proof of Lemma \ref{TaylorExpansion}, 
we see $f(-t)$ and  $g(t)$ can be given as follows.

\begin{corol}\label{Taylor-fg(t)}
\begin{align}
f(-t)&= \sum_{k \geq 0} \frac {(-1)^k}{k!}
B_+\lp \lb H_t(z) \pz \rb^k \rp, \label{Taylor-fg(t)-e1} \\
g(t)&= \sum_{k \geq 0} \frac {1}{k!}
B_+\lp \lb M_t(z) \pz \rb^k \rp. \label{Taylor-fg(t)-e2}
\end{align}
\end{corol}

\begin{lemma}\label{L2.2.2}
\begin{align}
g(t)f(-t)&=f(-t)g(t)=1. \label{DUE-1}
\end{align}
\end{lemma}
\pf For any $u(z)\in \kzz$, by 
Eqs.\,(\ref{NewTaylorExpansion-e1}), (\ref{NewTaylorExpansion-e2}) and the fact 
that $G_t(z)=F_t^{-1}(z)$, we have
\begin{align*}
g(t)f(-t)u(z)=g(t)u(F_t(z))=u(F_t(G_t(z)))=u(z).
\end{align*}
Hence, as elements of $\cDazz[[t]]$, we have $g(t)f(-t)=id=1$. 
Similarly, we can show $f(-t)g(t)=1$.
\epfv

Next, let us consider the 
{\it D-Log} of $F_t(z)$, which 
has been studied in \cite{E1}--\cite{E3}, \cite{N}, 
\cite{Z-exp} and \cite{WZ} for the commutative case.

\begin{lemma}\label{L2.2.3}
For any $\alpha \geq 1$, $F_t(z)\in \ataz$, 
there exists a unique 
$a_t(z)\in \kttzz^{\times n}$ 
with $a_{t=0}(z)=0$ and $o(a_t(z))\geq \alpha$
such that
\begin{align}\label{L2.2.3-e1}
e^{\lb a_t(z)\pz \rb}\, \cdot z = F_t(z),
\end{align}
where, as usual, the exponential in the equation 
above is given by 
\begin{align}\label{L2.2.3-e2}
e^{\lb a_t(z)\pz \rb} = \sum_{m\geq 0} \frac {1}{m!}\lb a_t(z)\pz \rb^m.
\end{align}
\end{lemma}

For the proof of the lemma, 
we refer the reader to
the proof of Proposition $2.1$ 
in \cite{Z-exp} which gives 
an elementary proof of this result 
for the commutative case.
The main idea of the proof 
is to solve the homogeneous (in $z$) parts of 
$a_t(z)$ recursively from Eq.\,$(\ref{L2.2.3-e1})$. 
Even though, the proof in \cite{Z-exp} is for the  
commutative case over $\bC$, it works equally well 
for the noncommutative case as long as 
the base algebra, which is $K[[t]]$ in our case, 
is a $\bQ$-algebra.

Following the terminology used in 
\cite{WZ} for the commutative case, 
we call $a_t(z)$ the {\it D-Log} of $F_t(z)$. Now, from 
the D-Log $a_t(z)$, we define a sequence $\{\phi_m \in 
\cD er^{[\alpha]}\langle \langle z \rangle \rangle \, | m \geq 1\}$ 
of the derivations of $\kzz$ by requiring 
 
\begin{align}\label{Def-d(t)}
d(t):=\sum_{m=1}^\infty \frac {t^{m}}m \phi_m = -\lb a_t(z)\pz \rb. 
\end{align}

\begin{lemma}\label{L2.2.4}
\begin{align}
e^{d(t)} = g(t). \label{DUE-2}
\end{align}
\end{lemma}

\pf First, note that, it is well-known that 
the exponential of any derivation 
of an algebra, when it is well-defined, 
is an automorphism of the algebra.
By this fact and 
Eqs.\,(\ref{L2.2.3-e1}) and (\ref{NewTaylorExpansion-e1}), 
we have, for any polynomial $u(z)$ in the free variables $z$,   
\begin{align*}
e^{\lb a_t(z)\pz \rb}u(z)&=u( e^{\lb a_t(z)\pz \rb}\, z) \\ 
&= u (F_t(z))\\
&=f(-t) u(z).
\end{align*}

Hence, as elements of $\cDazz[[t]]$, we must have 
\begin{align}
e^{\lb a_t(z)\pz \rb}=f(-t). \label{DUE-2-pe2}
\end{align}

Secondly, by Eq.\,(\ref{Def-d(t)}), 
we have $e^{\lb a_t(z)\pz \rb}=e^{-d(t)}$ 
whose inverse element in $\cDzz [[t]]$ is obviously $e^{d(t)}$; 
while by Eq.\,(\ref{DUE-1}), the inverse element of 
$f(-t)$ is given by $g(t)$.
Hence, by Eq.\,(\ref{DUE-2-pe2}) above, we have  $e^{d(t)}=g(t)$.
\epfv

One remark on the D-Log is as follows. 
Note that, by taking the D-Log, 
we have a well-defined map 
\begin{align}\label{T-map}
T: \ataz & \to \, \, t \, \cDrtazz  \\
F_t(z)& \to \quad  \lb a_t(z)\pz \rb. \nno
\end{align}

Conversely, for any $\delta_t \in t\,\cDrtazz =t \, \cDrazz [[t]]$, 
we can always write $\delta_t$ 
uniquely as $\delta_t =\lb a_t(z)\pz \rb$ for 
some $a_t(z)\in \kttzz$ with $a_{t=0}(z)=0$ and 
$o(a_t(z)) \geq \alpha$. 
By taking the exponential 
$e^{\delta_t} =e^{\lb a_t(z)\pz \rb}$ of 
the derivation $\delta_t$, 
we get an automorphism of $\kttzz$. 
From Eq.\,(\ref{L2.2.3-e2}), 
it is easy to see that, 
the resulting  automorphism 
$e^{\lb a_t(z)\pz \rb}$ actually lies in $\ataz$. 
Therefore we have the following proposition.

\begin{propo}\label{T-bijection}
For any $\alpha \geq 1$, 
$T: \ataz  \to  t \, \cDrtazz$ 
defined in Eq.\,$(\ref{T-map})$ is a bijection, 
whose inverse map is given by the exponential map.
\end{propo}

By using the Baker-Campbell-Hausdorff formula 
(see \cite{Va} or \cite{Re}), 
it is easy to check that, formally,
the Lie algebra of the Lie group $\ataz$ 
is exactly the Lie algebra $t\cDrtazz$.
The inverse map of $T$, 
which is the exponential map 
given by Eq.\,(\ref{L2.2.3-e2})
is same as the exponential map of 
the Lie group $\ataz$ 
from its Lie algebra 
to itself. So, in terms of the language 
of Lie theory, the proposition 
above just says that the exponential 
map of the formal Lie group $\ataz$ 
turns out to be a bijection 
whose inverse map is given by taking 
the D-Log of the elements of $\ataz$.

Next, we define the last two 
sequences $\{ \psi_m\,|\, m\geq 1\}$ 
and $\{ \xi_m\,|\, m\geq 1\}$ 
of elements of $\cDrazz$ by requiring 
\begin{align}
h(t)&:=\sum_{m\geq 1} \psi_m t^{m-1}=
\left[\frac{\p M_t}{\p t}(F_t) \pz \right ],\label{Def-h(t)}\\
m(t)&:=\sum_{m\geq 1} \xi_m t^{m-1}=
\left[\frac{\p H_t}{\p t}(G_t) \pz \right ]. \label{Def-m(t)}
\end{align}

With the facts that $o( H_t(z) )\geq \alpha$ and 
$o( M_t(z) )\geq \alpha$,
It is easy to see that the derivations 
$\psi_m$ and $\xi_m$ $(m\geq 1)$ defined above 
are indeed in $\cDrazz$.

To get some concrete ideas for the differential operators 
defined in Eqs.\,$(\ref{Def-h(t)})$ and $(\ref{Def-m(t)})$, 
let us recall the following result proved in 
Lemma $4.1$ in \cite{NC-IVP} for 
the special automorphism $F_t\in \ataz$ 
of the form $F_t(z)=z-tH(z)$ with 
$H(z)$ independent on $t$, i.e. 
$H(z) \in \kzz^{\times n}$.

\begin{lemma} \label{NCIVP-L4.1.1}
For any $F_t \in \ataz$ of the form $F_t(z)=z-tH(z)$ 
as above, let $N_t(z)=t^{-1}M_t(z)$. 
Then we have

\begin{align}
m(t)&=\lb N_t(z)\pz \rb, \label{NCIVP-L4.1.1-e1} \\
h(t)&=\sum_{m\geq 1}  t^{m-1} \lb C_m(z)\pz \rb, \label{NCIVP-L4.1.1-e2}
\end{align}
where $C_m(z)\in \kzz^{\times n}$ $(m\geq 1)$ 
are defined recurrently by
\begin{align}
C_1(z)& =H(z),\label{NCIVP-L4.1.1-e3}\\
C_m(z)&=\lb C_{m-1}(z)\pz \rb H(z),\label{NCIVP-L4.1.1-e4}
\end{align}
for any $m\geq 2$.

Consequently, for any $m\geq 1$, the derivations 
$\psi_m$ and $\xi_m$ defined in 
Eqs.\,$(\ref{Def-m(t)})$ and $(\ref{Def-h(t)})$
are given by

\begin{align}
\psi_m &=\lb C_m(z)\pz \rb, \label{C4.1.2-e1}\\
\xi_m&=\lb N_{[m]}(z)\pz \rb, \label{C4.1.2-e2}
\end{align}
where $N_{[m]}(z)\in \kzz^{\times n}$ $(m\geq 1)$ 
is the coefficient of $t^{m-1}$ of $N_t(z)$.
\end{lemma}

By the mathematical induction on $m\geq 1$, 
it is easy to show that, 
when $z$ are commutative variables, 
we further have
\begin{align}
C_m(z)=(JH)^{m-1}H(z) \label{Special-Cm}
\end{align}
for any $m\geq 1$, where $JH$ is the Jacobian matrix of 
$H(z)\in K[[z]]^{\times n}$.

One remark about Lemma \ref{NCIVP-L4.1.1}
is as follows. For the automorphisms 
$F_t\in \ataz$ of the form $F_t(z)=z-tH(z)$ 
as above, the differential operators 
$\psi_m$ $(m\geq 1)$ have the simplest 
form. In particular, 
in the commutative case,
by Eq,\,(\ref{Special-Cm}),
they capture the nilpotence of the Jacobian 
matrix $JH$ in a nice way, namely, 
$JH$ is nilpotent iff $\psi_m=0$ for any $m\geq n$.
On the other hand, the differential 
operators $\xi_m$ $(m\geq 1)$ capture the inverse map 
$G_t(z)$ of $F_t(z)$ directly. 
Together with the main result 
of this section (see Theorem \ref{S-Correspondence} below), 
this leads to an interesting connection 
of the NCSF theory with the inversion problem 
(\cite{BCW}, \cite{E}), 
i.e. the problem that studies 
various properties of the inverse map $G_t$ from $F_t$.
This connection will be discussed more in 
Remark \ref{R3.2.13} at the end of this subsection.

Now let us consider the relations of 
the differential operators 
$h(t)$ and $m(t)$ defined in 
Eqs.\,$(\ref{Def-h(t)})$ and $(\ref{Def-m(t)})$, 
respectively, with the differential operator 
$g(t)$ defined by Eq.\,$(\ref{Def-g(t)})$
or (\ref{TaylorExpansion-e2}).

\begin{lemma}\label{DUE-3-4}
\begin{align}
\frac {d g(t)} {d t} &= g(t) h(t), \label{DUE-3}  \\
\frac {d g(t)} {d t} &=m(t) g(t). \label{DUE-4}
\end{align}
\end{lemma}
\pf First, by Proposition $3.3$ in \cite{NC-IVP},  
for any $u(z)\in \kzz$, the following equations hold.
\begin{align}
\frac {\p \, u(G_t )}{\p t} &=(h(t)u)(G_t), \label{DUE-3-4-pe1} \\
\frac {\p \, u(G_t )}{\p t} &=m(t)\,\, u(G_t), \label{DUE-3-4-pe2}
\end{align}
Now, by Eq.\,(\ref{NewTaylorExpansion-e2}), we can rewrite 
Eq.\,(\ref{DUE-3-4-pe1}) as follows.
\begin{align}
\frac {\p }{\p t} (g(t)u(z))&=  g(t)(h(t)u)(z), \nno \\
\frac {d g(t)}{d\, t} \, u(z) &= (g(t)h(t))\, u(z).\label{DUE-3-4-pe3}
\end{align}
Since Eq.\,(\ref{DUE-3-4-pe3}) above holds for any $u(z)\in \kzz$, 
$\frac {d g(t)} {d t}$ and $g(t) h(t)$, as elements of $\cDzz[[t]]$, 
must be same. Hence we have Eq.\,(\ref{DUE-3}).

To show Eq.\,(\ref{DUE-4}), again by Eq.\,(\ref{NewTaylorExpansion-e2}), 
we can rewrite Eq.\,(\ref{DUE-3-4-pe2}) as follows.
\begin{align*}
\frac {\p }{\p t} (g(t)u(z))&=  m(t)(g(t)u)(z), \nno \\
 \frac {d g(t)}{d\, t} \, u(z) &= (m(t)g(t))\, u(z).
\end{align*}

By the same reason as above,  $\frac {d g(t)} {d t}$ 
and $m(t) g(t)$ must be same.
Hence we get Eq.\,(\ref{DUE-4}).
\epfv

Now we can formulate the main result of this subsection. 
First, we set
\begin{align}\label{Def-Omega-Ft}
\Omega_{F_t}:=(f(t), \, g(t),\, d(t),\, h(t),\, m(t)) \in \cDazz[[t]]^{\times 5}.
\end{align}
Then, comparing Eqs.\,(\ref{DUE-1}), (\ref{DUE-2}), 
(\ref{DUE-3}) and (\ref{DUE-4}) with 
Eqs.\,(\ref{UE-1})--(\ref{UE-4}), respectively, 
and noting that, 
by Lemma \ref{TaylorExpansion},  
Eq.\,(\ref{UE-0}) is also satisfied by $f(t)$, 
we get the following theorem.

\begin{theo} \label{GTS-II-Main}
For any $\alpha \geq 1$ and $F_t(z)\in \ataz$, 
the $5$-tuple $\Omega_{F_t}$ defined 
in Eq.\,$(\ref{Def-Omega-Ft})$
forms a \cNcs system over 
the differential operator $K$-algebra $\cDazz$.
\end{theo}

Furthermore, note that, 
the coefficients $\psi_m$ $(m\geq 1)$ 
are all $K$-derivations, hence are 
primitive elements of
the Hopf algebra $\cDazz$ 
(see Eq.\,(\ref{Coprd-delta})). 
Then, by Theorem \ref{Universal}, we have 
the following correspondence 
between the NCSFs in the universal \cNcs system 
$(\cNsf, \Pi)$ and the differential operators 
defined in this subsection.

\begin{theo}\label{S-Correspondence}
For any  $\alpha \geq 1$ and $F_t \in \ataz$ 
as fixed before, there exists a unique 
$K$-Hopf algebra homomorphism 
$\cS_{F_t}: \cNsf \to \cDazz$  such that 
$\cS_{F_t}^{\times 5}(\Pi)=\Omega_{F_t}$. 

More precisely, we have the following 
correspondence between the
NCSFs in $\Pi$ and 
the differential operators 
in $\Omega_{F_t}$.  Namely, 
for any $m\geq 1$, we have
\begin{align}
\mathcal S_{F_t} (\Lambda_m) & =\lambda_m, \label{L-l} \\
\mathcal S_{F_t} (S_m) & =s_m, \label{S-s} \\
\mathcal S_{F_t} (\Psi_m) & =\psi_m, \label{Psi-psi}  \\
\mathcal S_{F_t} (\Phi_m)& =\phi_m, \label{Phi-phi}  \\
\mathcal S_{F_t} (\Xi_m)& =\xi_m. \label{Xi-xi} 
\end{align}
\end{theo}

Note that, by Proposition \ref{bases}, 
any one of Eqs.\,(\ref{L-l})--(\ref{Xi-xi}) 
in turn completely 
determines the homomorphism $\cS_{F_t}$ 
in the theorem above. 

\vskip2mm

One direct consequence of 
Theorem \ref{S-Correspondence} 
above is that, for each $F_t \in \ataz$, 
we get a so-called {\it specialization} $\cS_{F_t}$
(following the terminology 
in the theory of symmetric functions) of NCSFs 
by differential operators. From now on, 
we will call $\cS_{F_t}$ 
the {\it differential operator specialization} 
of NCSFs associated with the automorphism $F_t\in \ataz$. 

Finally, let us end this section with the following remark on
a direct consequence of Theorem \ref{S-Correspondence} to 
the study of the inversion problem (\cite{BCW}, \cite{E}). 

\begin{rmk}\label{R3.2.13}
By Lemma $\ref{NCIVP-L4.1.1}$ 
and the comments after Eq.\,$(\ref{Special-Cm})$, 
when $F_t$ has the form $F_t(z)=z-tH(z)$ with 
$H(z)\in \kzz^{\times n}$, the differential 
operators $\psi_m$'s and $\xi_m$'s 
 $($see Eqs.\,$(\ref{C4.1.2-e1})$, 
$(\ref{C4.1.2-e2})$
and $(\ref{Special-Cm}))$ 
becomes important for the study 
of the inversion problem 
$($\cite{BCW},\,\cite{E}$)$. 
On the other hand, as we pointed out earlier 
in Remark \ref{Motivation},
by applying the homomorphism $\cS_{F_t}$ 
and the correspondence in 
Theorem $\ref{S-Correspondence}$, 
we can transform the polynomial identities 
among the NCSFs $($\cite{G-T}, \cite{GTS-III}$)$
in the system $\Pi$ into polynomial
identities among the corresponding
differential operators in the system 
$\Omega_{F_t}$. Some of these identities  
may be used to study certain properties 
of the inverse map, the D-Log and the formal flow 
of $F_t(z)$. For more detailed study 
in this direction, see the followed paper 
\cite{GTS-III}. 
\end{rmk}

\renewcommand{\theequation}{\thesection.\arabic{equation}}
\renewcommand{\therema}{\thesection.\arabic{rema}}
\setcounter{equation}{0}
\setcounter{rema}{0}

\section{\bf Differential Operator Specializations of 
Noncommutative Symmetric Functions} \label{S4}

Let $\cDrazz$, $\cDazz$, $\ataz$, $F_t\in \ataz$ and 
$\cS_{F_t}$ as in the previous section. 
In this section, we study more properties 
and consequences of the differential 
operator specializations 
$\cS_{F_t}: \cNsf \to \cDazz$ $(F_t\in \ataz)$ 
given in Theorem \ref{S-Correspondence}. 
First, we show in Proposition \ref{cS-graded}
that, for any $F_t\in \ataz$, 
the specialization $\cSft$ is
a homomorphism of graded $K$-Hopf algebras from $\cNsf$ 
to the subalgebra $\cDaz:=\cDz \cap \cDazz$, 
where $\cDz$ is the differential operator algebra
of the polynomial algebra $\kz$,  if and only if 
$F_t(z)=t^{-1}F(tz)$ for some automorphism $F(z)$ of $\kzz$.
Consequently, for any $F_t$ satisfying 
the condition above,  by taking the graded duals, 
we get a graded $K$-Hopf algebra homomorphism
$\cSft^*: \cDaz^*  \to \cQf$ from the
graded dual $\cDaz^*$ 
of $\cDaz$ to the Hopf algebra $\cQf$ 
(\cite{Ge}, \cite{MR} and \cite{St2})
of quasi-symmetric functions 
(see Corollary \ref{S-cQsf}). 
Secondly, we show in Theorem \ref{Isomorphism} 
that, with a properly defined 
group product $\circledast$ on the set 
$\text{\bf Hopf}\,(\cNsf, \cDazz)$
of all $K$-Hopf algebra homomorphisms 
from $\cNsf$ to $\cDazz$, the correspondence of 
$F_t \in \ataz$ and 
$\cS_{F_t} \in \text{\bf Hopf}\,(\cNsf, \cDazz)$
gives an isomorphism of groups. 
Finally, we show in Theorem \ref{StabInjc}
that the family of 
the specializations $\cS_{F_t}$ 
(with all $n\geq 1$ and $F_t \in \ataz$)
can distinguish any two different NCSFs. 

First let us fix the following notations.

Let $\cDz$ be the differential operator 
algebra of the polynomial
algebra $\kz$, i.e. $\cDz$ is the unital subalgebra of 
$\text{End}_K(\kz)$ generated 
by all $K$-derivations of $\kz$.  
For any $m\geq 0$, let $\cD_{[m]}\langle z \rangle$ 
be the set of all differential operators 
$U$ such that, for any homogeneous polynomial 
$h(z)\in \kz$ of degree $d\geq 0$, $U h(z)$ 
either is zero or is homogeneous of degree $m+d$. 
For any $\alpha\geq 1$, set $\cDaz:=\cDz\cap \cDazz$.
Then, we have the grading
\begin{align}\label{Grading-cDz}
\cDaz &=\bigoplus_{m \geq \alpha-1} 
\cD_{[m]} \langle z \rangle, 
\end{align}
with respect to which $\cDaz$ 
becomes a graded $K$-Hopf algebra.

Now, for any $\alpha \geq 2$, we let 
$\gtaz$ be the set of all automorphisms 
$F_t\in \ataz$ such that $F_t(z)=t^{-1}F(tz)$
for some automorphism $F(z)$ of $\kzz$.
It is easy to check that $\gtaz$ 
is a subgroup of $\ataz$.

\begin{propo}\label{cS-graded}
For any $\alpha\geq 2$ and $F_t\in \ataz$, 
the differential operator specialization 
$\cS_{F_t}$ is a graded $K$-Hopf algebra 
homomorphism $\cS_{F_t}: \cNsf\to \cDaz\subset \cDazz$ 
iff $F_t\in \gtaz$.
\end{propo}

\pf First, by the definition of $\gtaz$, it is easy to see that
$F_t(z)=z-H_t(z) \in \gtaz$ iff
$H_t(z)$ can be written as
\begin{align}\label{cS-graded-pe1}
H_t(z)= \sum_{m\geq 1} t^m H_{[m]}(z)
\end{align}
with $H_{[m]}(z)$ $(m\geq 1)$ homogeneous of degree $m+1$.
Hence, we have $F_t=z-H_t(z) \in \gtaz$ iff
the differential operator $\lb H_t(z) \pz\rb$ 
can be written as  
\begin{align}\label{cS-graded-pe2}
\lb H_t(z) \pz\rb = \sum_{m\geq 1} t^m \delta_m 
\end{align}
with $\delta_m=\lb H_{[m]}(z) \pz\rb \in \cD_{[m]}\langle z \rangle$ 
for all $m\geq 1$.

Now, we first assume that $\cS_{F_t}$ in Theorem 
$\ref{S-Correspondence}$ 
is a graded $K$-Hopf algebra homomorphism from 
$\cNsf$ to $\cD \langle z \rangle$, then, combining with Eq.\,(\ref{L-l}),
we have $\cSft(\Lambda_m)=\lambda_m \in \cD_{[m]}\langle z \rangle$ 
for any $m\geq 1$. But,  by Eq.\,(\ref{NewTaylorExpansion-e1}) 
with $u(z)=z$, we have $F_t(z)=f(-t)\cdot z$, i.e. 
$H_t(z)=\sum_{m\geq 1}(-1)^{m-1} t^m \lambda_m z$ with 
$\lambda_m z$ homogeneous of degree $m+1$. Therefore
Eq.\,(\ref{cS-graded-pe1}) holds with $H_{[m]}(z)=(-1)^{m-1}\lambda_m z$ 
homogeneous of degree $m+1$ for any $m\geq 1$. 
Hence, we have $F_t\in \gtaz$. 

Conversely, assume $F_t(z)\in \gtaz$. Then Eq.\,(\ref{cS-graded-pe2}) 
holds with $\delta_m\in \cD_{[m]}\langle z \rangle$ 
for all $m\geq 1$. Then by Eq.\,(\ref{Taylor-fg(t)-e1}),
it is easy to see that $\lambda_m \in \cD_{[m]} \langle z \rangle$ 
for all $m\geq 1$. By Eq.\,(\ref{L-l}) and
the fact that $\cNsf$ is the free $K$-algebra 
generated by $\Lambda_m$ $(m\geq 1)$, 
it is easy to see that $\cS_{F_t}$ does preserve
the gradings of $\cNsf$ and $\cDz$ defined 
by Eqs.\,$(\ref{Grading-cNsf})$ 
and $(\ref{Grading-cDz})$, respectively.   
\epfv

Now, for any $F_t\in \gtaz$ $(\alpha\geq 2)$, 
by the proposition above, 
we can take the graded dual of 
the graded $K$-Hopf algebra homomorphism 
$\cSft:\cNsf \to \cDaz$ 
and get the following corollary.

\begin{corol}\label{S-cQsf}
For any $\alpha \geq 2$ and $F_t\in\gtaz$, 
let $\cDaz^*$ be the graded dual 
of the graded $K$-Hopf algebra $\cDaz$. Then, 
$$
\cSft^*: \cDaz^* \to \cQf
$$ 
is a homomorphism of 
graded $K$-Hopf algebras. 
\end{corol}

Next, let us consider the following observation.
Note that, for any $F_t\in \ataz$, by looking at the 
differential operator specialization $\cS_{F_t}$
in Theorem \ref{S-Correspondence}, we get
the following map.

\begin{align}\label{Embedding}
\mathbb S: \ataz & \longrightarrow \text{\bf Hopf}\,(\cNsf, \cDazz)  \\
F_t \quad  & \longrightarrow \quad\quad\quad  \cS_{F_t}, \nno 
\end{align}
where $\text{\bf Hopf}\,(\cNsf, \cDazz)$ denotes 
the set of $K$-Hopf   algebra 
homomorphisms from $\cNsf$ to $\cDazz$. 

We claim the map $\mathbb S$ above is a bijection. 
To see it is injective, 
first note that $\cS_{F_t}$ is uniquely determined 
by $f(t)=\cS_{F_t}(\lambda(t))$ since
$\cNsf$ is freely generated 
by the coefficients $\Lambda_m$ 
$(m\geq 1)$ of $\lambda(t)$.
While by Eq.\,(\ref{NewTaylorExpansion-e1}) with $u(z)=z$, 
we have $F_t(z)=f(-t)\cdot z$. Therefore,
$f(t)$ conversely completely determines 
$F_t$ itself. Hence $\mathbb S$ is injective.
To show the surjectivity of $\mathbb S$,
let $\cT$ be any element of 
$\text{\bf Hopf}\,(\cNsf, \cDazz)$.
Then $\cT( \Phi (t))$ must have all 
its coefficients primitive in $\cDazz$. 
This is because the NCSFs 
$\Phi_m$ $(m\geq 1)$ are all primitive and any 
Hopf algebra homomorphism preserves 
primitive elements. 
Since the only primitive elements 
of $\cDazz$ are $K$-derivations,
we have $\cT(\Phi(t))\in t \,\cDrazz[[t]]=t\,\cDrtazz$.
Then by Proposition \ref{T-bijection}, 
there exist a unique
$F_t(z)\in \ataz$ such that its  
D-Log $a_t(z)$ is given by the property
$\lb a_t(z) \pz \rb=-\cT(\Phi(t))$.  
By Eq.\,(\ref{Phi-phi}) in 
Theorem \ref{S-Correspondence},
we have
$\cS_{F_t}(\Phi(t))=\cT(\Phi(t))$, 
hence, we must have $\cS_{F_t}=\cT$ since, by Proposition \ref{bases}, 
$\cNsf$ is also freely generated by $\Phi_m$ $(m\geq 1)$. 
So $\mathbb S$ is also surjective.

Next, we consider the following group product 
on the set $\text{\bf Hopf}$ $(\cNsf, \cDazz)$,
with which the bijection 
$\mathbb S$ in Eq.\,(\ref{Embedding}) 
becomes an isomorphism of groups. 
Note that the convolution product $\ast$ 
of the linear maps from $\cNsf$ 
to $\cDazz$ is not a right choice, 
for $\text{\bf Hopf}\,(\cNsf, \cDazz)$
is not closed under 
the convolution product.

First, let ${\bf DPS}\,(\cDazz)$ be the set 
of all sequences of divided powers 
$\{a_m\,|\, m\geq 0\}$ (with $a_0=1$) of the Hopf 
algebra $\cDazz$. 
Consider the following map:
\begin{align}\label{Embedding-2}
\mathbb D: \text{\bf Hopf}\,(\cNsf, \cDazz)  \rightarrow & \quad \, {\bf DPS}\,(\cDazz) \\
\cT \quad \quad \quad \longrightarrow &\quad \{ \cT (S_m)\,|\, m\geq 0 \}. \nno 
\end{align}

Note that the map $\mathbb D$ is well-defined 
for the NCSFs $\{S_m\,|\, m\geq 0\}$ form 
a sequence of divided powers of $\cNsf$ 
and any Hopf algebra homomorphism 
preserves sequences of divided powers. 

We claim the map $\mathbb D$ 
above is a bijection. 
The injectivity is obvious for $\cNsf$ 
is the free $K$-algebra 
generated by $S_m$ $(m\geq 1)$ 
(see Proposition \ref{bases}).
To see the surjectivity of $\mathbb D$, 
let $\{a_m\,|\, m\geq 0\}\in {\bf DPS}\,(\cDazz)$ 
and $\mathcal A: \cNsf\to \cDazz$ 
the unique $K$-algebra 
homomorphism 
that maps $S_m$ to $a_m$ 
for any $m\geq 1$. 
Then, by applying Theorem 
\ref{Universal}, $(b)$ 
to the \cNcs system 
$\mathcal A^{\times 5}(\Pi)$ over $\cDazz$, 
we see that $\mathcal A$ must also be a $K$-Hopf 
algebra homomorphism, i.e.
$\mathcal A\in \text{\bf Hopf}\,(\cNsf, \cDazz)$.
By the definition of $\mathcal A$, 
we have $\mathbb D(\mathcal A)=\{a_m\,|\, m\geq 0\}$. 
Hence $\mathbb D$ is also surjective.

Now we identify any sequence of divided powers  
$\{a_m\,|\, m\geq 0\}$ of $\cDazz$ 
with its generating function 
$a(t):=\sum_{m\geq 0} a_m t^m$. Note that, 
in general, a sequence $\{a_m\,|\, m\geq 0\}$
with $(a_0=1)$ is a divide series 
iff its generating function $a(t)$ satisfies 
$\Delta a(t)=a(t)\otimes a(t)$ 
(and $a_0=1$). 
Therefore we can view ${\bf DPS}\,(\cDazz)$ as 
the subset of the elements 
$a(t)$ of $\cDazz[[t]]$ such that
$\Delta a(t)=a(t)\otimes a(t)$ 
(and $a_0=1$). 
Note that, for any $a(t), b(t) \in {\bf DPS}\, 
(\cDazz)\subset \cDazz[[t]]$, we have
\begin{align*}
\Delta (a(t)b(t))&=\Delta (a(t)) \Delta(b(t))\\
&=(a(t)\otimes a(t))(b(t)\otimes b(t))\\
&=(a(t)b(t))\otimes (a(t)b(t))
\end{align*}
So, ${\bf DPS}\,(\cDazz)$ is closed under 
the algebra product of $\cDazz[[t]]$.
Since each $a(t)\in {\bf DPS}\,\,(\cDazz)$ 
is also invertible in $\cDazz[[t]]$,
for $a(0)=1$, 
${\bf DPS}\,(\cDazz)$ 
with the algebra product of $\cDazz[[t]]$
forms a group.
By identifying the set 
$\text{\bf Hopf}\,(\cNsf, \cDazz)$ with 
${\bf DPS}\,(\cDazz)$ via
the bijection $\mathbb D$ in Eq.\,(\ref{Embedding-2}),
we get a group product denoted by $\circledast$ for 
the set $\text{\bf Hopf}\,(\cNsf, \cDazz)$.
Let $\circ$ denote the group product of $\ataz$  
from the composition of automorphisms. 
Then we have the following theorem.

\begin{theo}\label{Isomorphism}
For any $\alpha \geq 1$, the map 
\begin{align}\label{Isomorphism-e1}
\mathbb S: (\, \ataz, \, \circ\,& )  
 \to (\, \text{\bf Hopf}\,(\cNsf, \cDazz), \, \circledast \, ) \\
F_t \quad \quad & \longrightarrow  \quad \quad \quad\quad \cS_{F_t}, \nno 
\end{align}
is an isomorphism of groups.
\end{theo}

\pf By the discussion after Eq.\,(\ref{Embedding}), 
we only need show  the map $\mathbb S$ 
is a homomorphism of groups. Furthermore, 
by the definition of the group product 
$\circledast$ above, it will be enough 
to show that, for any $U_t(z), V_t(z)\in \ataz$,
\begin{align}
\cS_{U_t\circ V_t}(\sigma (t))
&=\cS_{U_t}(\sigma (t)) \cS_{V_t}(\sigma (t)), \label{Isomorphism-pe1}\\
\cS_{U^{-1}_t}(\sigma (t))&=\cS_{U_t}(\sigma (t))^{-1}, \label{Isomorphism-pe2}
\end{align}
where, as before, $\sigma (t)=\sum_{m\geq 0} S_m t^m$ is 
the generating function of 
the complete elementary homogeneous NCSFs 
$\{ S_m \,|\, m\geq 0\}$.

First, by Theorem \ref{S-Correspondence},
we know that, for any $F_t(z)\in \ataz$, 
$\cS_{F_t}(\sigma (t))=g(t)$, which is the 
differential operator defined in Eq.\,(\ref{Def-g(t)}).  
Therefore,  
for any $v_t(z)\in \kttzz$, by Eq.\,(\ref{NewTaylorExpansion-e2}),
we have
\begin{align}
\cS_{F_t}(\sigma (t))v_t(z) = v_t (F^{-1}_t(z)). \label{Isomorphism-pe3}
\end{align}
In particular, for any $u(z) \in\kzz$, we have 
\begin{align*}
\cS_{U_t}(\sigma (t)) \cS_{V_t}(\sigma (t)) \, u(z)
&=\cS_{U_t}(\sigma (t)) \,u(V^{-1}_t(z))\\
&=u (\,V^{-1}_t \circ U^{-1}_t (z)\, )\\
&=u(\,(U_t\circ V_t)^{-1} (z)\,)\\
&=\cS_{U_t\circ V_t}(\sigma (t))\, u(z).
\end{align*}
Hence, we have Eq.\,(\ref{Isomorphism-pe1}).

To show Eq.\,(\ref{Isomorphism-pe2}), applying 
Eq.\,(\ref{Isomorphism-pe3}) to $F_t=U^{-1}_t$ and
$v_t(z)=u(z)$ for any $u(z)\in \kzz$, 
we have, 
\begin{align}
\cS_{U^{-1}_t}(\sigma (t))u (z)=u (U_t(z)). \label{Isomorphism-pe4}
\end{align}
On the other hand, by Theorem \ref{S-Correspondence} 
and Eq.\,(\ref{NewTaylorExpansion-e1}), we also have
\begin{align*}
\cS_{U_t}( \sigma(t) )^{-1}u(z)&= \cS_{U_t}( \sigma(t)^{-1} )u(z)\\
&=\cS_{U_t}(\lambda(-t)) u(z) \\
&=f(-t) u(z)\\
&=u(U_t(z)). 
\end{align*}
Hence, we have Eq.\,(\ref{Isomorphism-pe2}).
\epfv

\begin{rmk}\label{R3.2.16}
Note that, one implication of 
Theorem \ref{Isomorphism} 
and Proposition \ref{cS-graded} is that, 
all specializations 
$\cT:\cNsf \to \cDazz$ 
of NCSFs, which are also 
 $K$-Hopf algebra homomorphisms, 
are exactly the specializations 
$\cSft$ $(F_t\in \ataz)$;
while all specializations 
$\cT:\cNsf \to \cDaz$ 
of NCSFs, which are also graded
 $K$-Hopf algebra homomorphisms, 
are exactly the specializations 
$\cSft$ with $F_t\in \gtaz$.
\end{rmk}

Next we show that, for any $\alpha\geq 2$, 
the family of the differential operator 
specializations $\cS_{F_t}$ 
(with all $n\geq 1$ and 
$F_t\in \ataz$) can distinguish 
any two different NCSFs. 
This statement obviously is same as 
the following theorem. 

\begin{theo}\label{StabInjc}
In both commutative and noncommutative cases, 
the following statement holds. 
\vskip2mm
For any fixed $\alpha\geq 2$ and non-zero $P \in \cNsf$, 
there exist $n\geq 1$ 
$($the number  of the free variable 
$z_i$'s$)$ and $F_t(z)\in \ataz$ such that 
$\cS_{F_t} (P)\neq 0$.
\end{theo}
 
This is probably the only result 
derived in this paper for 
which the commutative case 
does not follow from the noncommutative 
case by the procedure of abelianization. 
Instead, as one can easily show that
the noncommutative case actually follows
from the commutative case by 
choosing any lifting of $F_t$ 
in commutative variables 
to an automorphism of formal 
power series algebra in 
noncommutative variables.
Therefore, we need only prove the theorem 
for the commutative case.

\vskip2mm

\pf First, by Proposition \ref{bases}, 
we may view $\cNsf$ as the free algebra $\kphi$
generated by $\Phi_m$ $(m\geq 1)$. 
Below we will write $\cNsf$ 
as $\kphi$ and view any NCSF 
$P\in \cNsf=\kphi$ as a polynomial 
$P(\Phi)$ in $\Phi_m$'s. 
Secondly, let $\mathcal L(\Phi) \subset \kphi$ be 
the free Lie algebra generated by 
$\Phi_m$ $(m\geq 1)$. 
Then, for any free variables $z$, 
by a similar argument as the proof of
the bijectivity of 
the map $\mathbb D$ in 
Eq.\,(\ref{Embedding-2}),
it is easy to see that, via the restriction 
and extension of the homomorphisms,
the set $\text{\bf Hopf}\, (\cNsf,  \cDazz )$ 
is in $1$-$1$ correspondence with 
the set $\text{\bf Lie}(\mathcal L(\Phi), \, \cDrazz )$
of the Lie algebra homomorphisms 
from $ \mathcal L(\Phi)$ to 
$\cDrazz$. 
Combining with the isomorphism 
$\mathbb S: \ataz \simeq 
\text{\bf Hopf}\, (\cNsf,  \cDrazz )$ 
in Theorem \ref{Isomorphism}, 
we see that the set of the specializations 
$\cS_{F_t}$ with $F_t\in \ataz $ is 
in $1$-$1$ correspondence with the set 
$\text{\bf Lie}(\mathcal L(\Phi), \, \cDrazz )$.

Now let $\mathcal K$ be the set of all 
NCSFs $P(\Phi)\in \kphi$ that
are mapped to zero by
the specialization $\cS_{F_t}$ 
for any number $n\geq 1$ of free 
variables $z=(z_1, z_2, ... , z_n)$ 
and any $F_t(z)\in \ataz$.
By the $1$-$1$ correspondence discussed above,
it is easy to see that $P(\Phi) \in \mathcal K$ 
iff $P(\Phi)$ satisfies the following property:

\vskip2mm

({\bf K})   {\it 
For any number $n\geq 1$ of free variables 
$z=(z_1, z_2, ... , z_n)$ and any 
Lie algebra homomorphism $\sigma: \mathcal L(\Phi) \to 
\cDrazz $, if
we still denote by $\sigma$
the extended $K$-algebra homomorphism  
from $\kphi$ to $\cDazz$. 
Then we have $\sigma (P(\Phi))=0$.}

\vskip2mm

Now we assume that the theorem is false, i.e. $\mathcal K \neq 0$ 
and derive a contradiction as follows. 

\vskip2mm

{\bf Claim $1$:} {\it $\mathcal K$ is invariant 
under linear transformations. More precisely, 
for any integers $M, N\geq 1$ and $a_{i, j}\in K$ 
with $1\leq i\leq M$ 
and $1\leq j\leq N$, let $Y:=\{Y_m \,|\, m\geq 1\}$, 
where
\begin{align}
Y_m:=\begin{cases} \sum_{j=1}^N a_{m, j} \Phi_j \quad & \text{if \, $1\leq m\leq M$}, \\
\Phi_m  \quad & \text{if \, $m>M$}.
\end{cases}
\end{align}
Then, for any $P(\Phi) \in \mathcal K$, we have
$P(Y)\in \mathcal K$.}

\vskip2mm

\underline{\it Proof of Claim $1$:} It will be enough to show that 
$P(Y)$ satisfies the property ({\bf K}) when $P(\Phi)$ does.  
For any 
$\sigma \in \text{\bf Lie}(\mathcal L(\Phi), \, \cDrazz )$,
we define a homomorphism $\eta: \mathcal L(\Phi)  \to   \cDrazz$
of Lie algebras 
by setting $\eta(\Phi_m)=\sigma(Y_m)$ 
for any $m\geq 1$. 
Then it is easy to see that $\sigma (P(Y))=\eta (P(\Phi))$. 
Hence, 
when $P(\Phi)$ satisfies the property ({\bf K}), 
we have $\sigma (P(Y))=\eta (P(\Phi))=0$ for any 
$\sigma \in \text{\bf Lie}(\mathcal L(\Phi), \, \cDrazz )$.
\epfv

{\bf Claim $2$:} {\it Let $\mathcal H$ be the set of all elements 
$Q(\Phi) \in \mathcal K$ such that $Q(\Phi)$ is homogeneous 
in each $\Phi_m$ that is involved in $Q(\Phi)$. 
Then $\mathcal H \neq 0$.}

\vskip2mm

\underline{\it Proof of Claim $2$:}
Let $P(\Phi)$ be any non-zero element of 
$\mathcal K$. By Claim $1$ above,
we may assume that $P(\Phi)$ 
involves exactly 
$\Phi_m$ $(1\leq m\leq N)$ 
for some $N\geq 1$. 
Let $y=(y_1, y_2, ... , y_N)$ be $N$
formal central parameters, i.e. 
they commute with each other, and also 
with $\Phi_m$'s and all free 
variables $z_i$'s 
under the consideration. 
Let $P(\Phi; y)$ be the polynomial 
in $\Phi$ and $y$ obtained 
by replacing $\Phi_m$ $(1\leq m\leq N)$ 
in $P(\Phi)$ by $y_m \Phi_m$. 
We view $P(\Phi; y)$ 
as a polynomial in $y$ 
with coefficients in 
$\kphi$ and write it as 
\begin{align}\label{StabInjc-pe1}
P(\Phi; y)=\sum_{I\in \bN^N}  P_I(\Phi) \, y^I. 
\end{align}

Now,  for any $\Vec{v} \in K^{\times N}$, 
by Claim $1$, we have $P(\Phi; \Vec{v}) \in \mathcal K$. 
Therefore, $\sigma (P(\Phi; \Vec{v}))=0$ for 
any $\sigma \in \text{\bf Lie}
(\mathcal L(\Phi), \, \cDrazz )$. 
In particular, for any 
$u(z) \in \kzz $,  
we have $\sigma (P(\Phi; \Vec{v} ))\cdot u(z)= 0$. 
Therefore, when we write
$\sigma (P(\Phi; y ))\cdot u(z)=\sum_{I\in \bN^N} 
y^I \sigma(P_I(\Phi))\cdot u(z)$ 
as a formal power series in $z$ 
with coefficients 
in $K[y]$, all its coefficients 
will vanish at any $\Vec{v} \in K^{\times N}$. 
Hence, as polynomials in 
the commutative variables $y$, 
these coefficients must be 
identically zero, since $K^N$ 
is obviously dense with respect to 
the Zariski topology of $K^N$. 
Therefore, for any $u(z)\in \kzz$,
$\sigma (P(\Phi; y )) u(z)=0$ 
as an element of 
$K[y]\langle\langle z\rangle\rangle$.
Hence, $\sigma (P(\Phi; y ))=0$ 
as a polynomial of $y$ with coefficients in $\cDazz$. 
In particular, as elements of $\cDazz$,
all the coefficients $\sigma(P_I(\Phi))$ 
of $y^I$ in $\sigma (P(\Phi; y ))$ 
are equal to zero. 
Since this is true 
for any $\sigma \in \text{\bf Lie}
(\mathcal L(\Phi), \, \cDrazz )$,
all the coefficients $P_I(\Phi)$ 
of $y^I$ in Eq.\,(\ref{StabInjc-pe1}) 
are in $\mathcal K$. 

On the other hand,  by the definition of 
$P(\Phi; y)$ above, it is easy to see that, 
for any fixed $I=(i_1, \dots , i_N) \in \bN^N$, 
$P_I(\Phi)$ in Eq.\,(\ref{StabInjc-pe1}) 
is homogeneous in each $\Phi_m$ $(1\leq m\leq N)$ 
of partial degree $i_m$. 
Therefore, all $P_I (\Phi) \in \mathcal H$. 
Since not all $P_I(\Phi)$ can be zero, 
otherwise $P(\Phi)$ 
would be zero, hence
Claim $2$ follows. \epfv

{\bf Claim $3$:} {\it Let $\mathcal H_1$ be the set of all elements 
$Q(\Phi) \in \mathcal K$ such that $Q(\Phi)$ is homogeneous 
of partial degree $1$ in each $\Phi_m$ 
that is involved in $Q(\Phi)$. 
Then $\mathcal H_1 \neq 0$.}

\vskip2mm

\underline{\it Proof of Claim $3$:}
By using Claim 2, we first fix a nonzero 
element $P(\Phi)\in \mathcal H$ and, 
by Claim 1, we assume that $P(\Phi)$ 
involves exactly
$\Phi_m$ $(1\leq m\leq N)$ 
for some $N\geq 1$. 
For each $1\leq m\leq N$, we assume that
$P(\Phi)$ is homogeneous of degree 
$d_m\geq 1$ in $\Phi_m$. Hence $P(\Phi)$ 
is homogeneous of total degree 
$d:=\sum_{1\leq m\leq N} d_m$.

Let $B=\{(m, j)|  1\leq m\leq N;   
1\leq j\leq d_m\}$ and
$W=\{ W_{m, j}  |   (m, j)\in B \}$ 
any subset of $\{\Phi_m\,|\, m\geq 1\}$ 
with $|W|=d$.
Let $y=\{ y_{m, j} \, | \,  (m, j)\in B \}$ 
be a family of central formal parameters. 
For any $1\leq m\leq N$, 
set $Y_m:=\sum_{j=1}^{d_m} y_{m, j} W_{m, j}$  
and $Y:=\{Y_m\,|\, 1\leq m\leq N \}$. 
Let $\tilde P(W, y):=P(Y)$.
We view $\tilde P(W, y)$ as 
a polynomial in $y$ with coefficients 
in $K\langle W \rangle$ and let $Q(W)$ 
be the coefficient of the monomial 
$\prod_{(m, j)\in B} y_{m, j}$ 
in $\tilde P(W, y)$. 
First, by a similar argument 
as in the proof of 
Claim $2$ above, 
we see that $Q(W)$ as well as 
other coefficients of 
the monomials of $y$ appearing in 
$\tilde P(W, y)$ are in $\mathcal K$.
Secondly, it is easy to see that 
$Q(W)$ is homogeneous of partial degree $1$ 
in each $W_{m, j}$ $((m, j)\in B)$. 
Thirdly, $Q(W)\neq 0$ since 
$P(\Phi)$ can be recovered from $Q(W)$
by replacing $W_{m, j}$ by 
$\Phi_m$ for any $(m, j)\in B$ and then 
dividing the multiplicity factor $d_1!d_2!\cdots d_m!$.
Therefore, $Q(W)$ is a non-zero element of 
$\mathcal H_1$ and Claim $3$ follows. \epfv

Finally, we derive a contradiction as follows.

\vskip2mm

Let $P(\Phi)$ be a non-zero element of $\mathcal H_1$ 
with the least total degree in $\Phi$ 
among of all non-zero elements of $\mathcal H_1$. 
By Claim 1,  we may assume that $P(\Phi)$ 
is homogeneous of partial degree $1$ 
in each $\Phi_m$ $(1\leq m\leq N)$ 
for some $N\geq 1$ and does not depend 
on any $\Phi_m$ with $m>N$.
Now we write $P(\Phi)$ as 
\begin{align}\label{StabInjc-pe3}
P(\Phi)=\sum_{ m=1}^N P_m(\Phi) \, \Phi_m.
\end{align}

First,  not all $P_m(\Phi)$'s above 
can be zero since $P(\Phi)\neq 0$.  
Without losing any generality, 
let us assume $P_N(\Phi)\neq 0$.
Secondly, $P_N(\Phi)$ is homogeneous of partial degree $1$ 
in each $\Phi_m$ $(1\leq m\leq N-1)$ and 
does not depend 
on any $\Phi_m$ with $m\geq N$. Thirdly, 
the total degree 
$\deg P_N (\Phi)=\deg P(\Phi)-1 < \deg P (\Phi)$.

Now, for any $n\geq 1$. Let $z=(z_1, z_2, ... , z_n)$
be $n$ free variables
and $w$ a free variables
that is independent with $z$.  
For any $u(z)\in \kzz$ with $o(u(z))\geq \alpha$
and any $\sigma \in 
\text{\bf Lie}(\mathcal L(\Phi), \, \cDrazz )$,  
let $\tilde \sigma: \mathcal L(\Phi) \to 
\cD er^{[\alpha]} \langle\langle x; w\rangle\rangle$ 
be the unique Lie algebra homomorphism  
which maps $\Phi_m$ $(1\leq m\leq N-1)$ to 
$\sigma(\Phi_m)$; $\Phi_N$ to $\lb u(z)\frac\p{\p w} \rb $ 
and $\Phi_m$ $(m>N)$ to zero.
Then, we have 
$\tilde \sigma ( P_N(\Phi)) = \sigma (P_N(\Phi))$ and, 
for  any $1\leq m \leq N-1$, $\tilde \sigma(\Phi_m )\cdot w
= \sigma( \Phi_m ) \cdot w=0$ for $\sigma(\Phi_m)\in \cDrazz$. 
By the fact $P(\Phi)\in \mathcal K$ and 
Eq.\,({\ref{StabInjc-pe3}),
we have
\begin{align}
0&=\tilde \sigma (P(\Phi))\cdot w \nno \\
&=\sigma( P_N(\Phi) ) \lb u(z)\frac\p{\p w}  \rb \cdot w \nno\\
&=\sigma( P_N(\Phi) ) u(z).\nno
\end{align} 
Since the differential operator 
$\sigma( P_N(\Phi) )$ annihilates any 
$u(z)\in \kzz$ with $o(u(z)) \geq \alpha$, 
by Lemma \ref{L3.1.5}, $\sigma( P_N(\Phi) )=0$.
Since this is true for any $\sigma \in 
\text{\bf Lie}(\mathcal L(\Phi), \, \cDrazz )$,  
hence, $P_N(\Phi)\in \mathcal K$. 
By the facts pointed 
in the previous paragraph,
we further have, $P_N(\Phi)\in \mathcal H_1$. 
But the total degree 
$\deg P_N(\Phi)< \deg P(\Phi)$,
which contradicts to the choice of $P(\Phi)$.
\epfv

Finally, let us end this paper with 
the following remarks. 

First, it is easy to see that,
for any $\alpha\geq 2$, 
all the results derived in this paper 
still hold if the base algebra $K[[t]]$ 
is replaced by $K[t]$. Secondly,
in the followed paper \cite{GTS-V}, 
by using Theorem \ref{StabInjc} above
and some connections of the \cNcs system 
$(\cDazz, \Oft)$ with 
the \cNcs systems constructed 
in \cite{GTS-V} over 
the Grossman-Larson 
Hopf algebra of labeled rooted trees, 
the following much stronger version 
of Theorem \ref{StabInjc} 
will be proved. 

Let $\mathbb B^{[\alpha]}_t\langle z \rangle$
be the set of automorphisms $F_t(z)=z-H_t(z)$ 
of the polynomial algebra 
$\ktz$ over $K[t]$ such that 
the following conditions are satisfied.
\begin{enumerate}
\item[$\bullet$]  $H_{t=0}(z)=0$. 
\item[$\bullet$] 
$H_t(z)$ is homogeneous in $z$ 
of degree $d \geq \alpha$.
\item[$\bullet$] 
With a proper permutation of the free variables $z_i$'s,
the Jacobian matrix $JH_t(z)$ is strictly lower triangular.
\end{enumerate}

\begin{theo}\label{StabInjc-best}$($\cite{GTS-V}$)$
In both commutative and noncommutative cases, 
the following statement holds. 
\vskip2mm
For any fixed $\alpha\geq 2$ and non-zero $P \in \cNsf$, 
there exist $n\geq 1$ 
$($the number  of the free variable 
$z_i$'s$)$ and $F_t(z)\in \mathbb B^{[\alpha]}_t\langle z\rangle$ 
such that $\cS_{F_t} (P)\neq 0$.
\end{theo}

{\small \sc Department of Mathematics, Illinois State University,
Normal, IL 61790-4520.}

{\em E-mail}: wzhao@ilstu.edu.

\end{document}